\documentclass[authoryear,12pt]{elsarticle}
\usepackage{cancel}
\usepackage{caption}
\usepackage{color}
\usepackage{lscape}
\usepackage{afterpage}
\usepackage{pstricks}
\usepackage{pst-plot}
\usepackage{longtable}
\usepackage{dcolumn}
\usepackage{pst-node}
\usepackage{amsmath}
\usepackage{amssymb}
\usepackage{amsthm}
\usepackage{multirow}
\usepackage{colortab}
\usepackage{color,soul}
\usepackage{array}
\usepackage{slashbox}
\usepackage{colortbl}
\usepackage{subfigure}
\usepackage{textcomp}
\usepackage{pstricks, pst-node}
\usepackage{pst-all}
\usepackage{algorithm,algorithmic}
\usepackage{url}
\usepackage{soul}
\usepackage{hyperref}
\psset{arrows=->, labelsep=3pt, mnode=circle}

\usepackage{tabularx}

\usepackage{eurosym}

\usepackage{tikz}

\newfont{\rams}{msbm10 scaled\magstep1}
\newcommand{\rea}{\mbox{\rams \symbol{'122}}}

\newenvironment{resumeT}{\begin{list}{}{\setlength{\rightmargin}{\leftmargin}}\item[]
{\centering {\bf \it~~~}
\par}\item[]\ignorespaces}{\unskip\end{list}}

\newtheorem{note}{Note}[section]

\begin{document}
\title{Inducing a probability distribution in Stochastic Multicriteria Acceptability Analysis}

\author[Eco]{\rm Sally Giuseppe Arcidiacono}
\ead{s.arcidiacono@unict.it}
\author[Eco]{\rm Salvatore Corrente}
\ead{salvatore.corrente@unict.it}
\author[Eco]{\rm Salvatore Greco}
\ead{salgreco@unict.it}

\address[Eco]{Department of Economics and Business, University of Catania, Corso Italia, 55, 95129  Catania, Italy}
%\address[por]{University of Portsmouth, Portsmouth Business School, Centre of Operations Research and Logistics (CORL), Richmond Building, Portland Street, Portsmouth PO1 3DE, United Kingdom}

\date{}
\maketitle

\vspace{-1cm}

\begin{resumeT}

\textbf{Abstract:} \noindent In multiple criteria decision aiding, very often the alternatives are compared by means of a value function compatible with the preferences expressed by the Decision Maker. The problem is that, in general, there is a plurality of compatible value functions, and providing a final recommendation on the problem at hand considering only one of them could be considered arbitrary to some extent. For such a reason, Stochastic Multicriteria Acceptability Analysis gives information in statistical terms by taking into account a sample of models compatible with the provided preferences. These statistics are given assuming the existence of a probability distribution in the space of value functions being defined a priori. In this paper, we propose some methods aiming to build a probability distribution on the space of value functions considering the preference information given by the Decision Maker. To prove the goodness of our proposal we performed an extensive set of simulations. Moreover, a sensitivity analysis on the variables of our procedure has been done as well. 

%\vspace{1cm}
%%%%%%%%%%%%%%%%%%%%
%%%%%%% Keywords %%%%%%%
%%%%%%%%%%%%%%%%%%%%
\end{resumeT}

%\vfill\newpage
\pagenumbering{arabic} 

%%%%%%%%%%%%%%%%%%%%%%%%%%%%%%%%%%%%
\section{Introduction}\label{intro}%
%%%%%%%%%%%%%%%%%%%%%%%%%%%%%%%%%%%%
In decision support problems we frequently face the issue of the lack of a complete information. This lack determines the tradeoff between the depth of the provided information and the accuracy of the recommendations. All decision models need the specification of some parameters to be implemented. Usually, such parameters strongly depend on the preferences of the decision maker (DM) and to determine them, one can adopt a direct or an indirect approach. The former requires that the necessary parameters are directly provided by the DM. So, this approach assumes that the DM clearly knows and understands the preference model and how the parameters affect the results. Unfortunately, in real-world decision problems, such an assumption is too demanding. On the contrary, in the indirect elicitation paradigm of the Ordinal Regression approach \citep{JacquetLagrezeSiskos1982}, we infer the values to be assigned to the parameters of the model such that it restores the preference information given by the DM. This approach is less demanding than the direct approach and it has been successfully used in literature as well as in several real-world applications {\citep{JacquetLagrezeSiskos2001,Kadzinski2022}. The simplicity of the indirect elicitation is however counterbalanced from the fact that in general, a plurality of models, e.g. a plurality of value functions, are compatible with the preferences provided by the DM and the application of each one of these models could give different results on the problem under consideration. For such a reason, the selection of one of these compatible models is arbitrary and meaningless to some extent. To overcome this choice problem, the Robust Ordinal Regression (ROR; \citealt{GrecoMousseauSlowinski2008,CorrenteEtAl2013}) and the Stochastic Multicriteria Acceptability Analysis (SMAA; \citealt{LahdelmaHokkanenSalminen1998,PelissariEtAl2020}) have been proposed in the literature. On the one hand, ROR takes into account simultaneously all the models compatible with the preferences given by the DM defining a necessary and a possible preference relation. Given two alternatives $a,b\in A$, where $A$ is the set of alternatives under consideration, $a$ is necessarily preferred to $b$ iff $a$ is at least as good as $b$ for all compatible models, while $a$ is possibly preferred to $b$ iff $a$ is at least as good as $b$ for at least one of them. On the other hand, the SMAA methodology gives information on the problem at hand in statistical terms taking into consideration some probability distribution in the space of compatible models (for example, in the space of the weight vectors defining the weighted sum preference model). From an operational point of view, a well-distributed sample of compatible models in the corresponding space is considered. To each of these models (for example, to each weighted sum characterized by a given weight vector) is assigned a certain ``probability" (in general equal for all of them), and, then, SMAA provides the rank acceptability indices for each alternative and each rank position and the pairwise winning index for each ordered pair of alternatives. The rank acceptability index gives the probability with which an alternative fills a certain specific rank position, while, the pairwise winning index gives the probability with which an alternative is preferred to another. 

As underlined above, in general, it is assumed that all sampled compatible models have the same probability. However, this is not necessarily true. In \cite{CorrenteEtAl2016KNOSYS} a procedure to estimate a probability distribution on the space of sampled compatible models has been proposed taking into account certain and uncertain preferences given by the DM. With certain preferences, we refer to preferences for which the DM has no doubt, such as ``I definitely prefer $a$ to $b$", while for uncertain preferences we refer to preferences for which the DM has a general tendency without expressing a definitive certitude, such as ``I have a certain tendency to prefer $a$ to $b$". In particular, in the spirit of the ordinal regression approach, the proposed procedure selects the probability distribution on the space of compatible value functions that maximizes the minimal difference between the probability of the preference of $a$ over $b$ on the probability of the preference of $b$ over $a$, for all pair of alternatives $(a,b)$ for which the DM has expressed a preference of $a$ over $b$.  Observe that differently from several other approaches considered in the literature, no assumption is made on the probability distribution in the space of value functions. In fact, a large variety of procedures have been presented to induce a probability distribution in the space of value functions on the basis of the Bayesian paradigm (see e.g. \citealt{chu2005preference,guo2010real, chajewska2013utilities,lepird2015bayesian,yet2020estimating,ru2022bayesian}). In this paper, we propose a different approach. Instead of inducing a probability distribution on the whole space of value functions, we induce a probability distribution on a well-distributed sample of value functions. In particular, on the basis of the preference information supplied by the DM in terms of pairwise comparisons of alternatives, a reference value function is determined and, maximizing the probability of the DM's preference comparisons,  to each value function in the sample is assigned a probability mass decreasing with the distance from the reference value function. This problem can be formulated in different forms permitting to use linear programming (in case the probability mass is assumed piecewise linear with respect to the distance from the reference function) or nonlinear optimization of single variable functions (in case the probability mass is assumed to be represented as some single parameter function of the distance from the reference value function). Therefore we use a rather simple procedure, obtaining in any case results that prove to be effective and reliable as we shall show in the computational experiments presented in the paper.   

%In this paper we are then interested in proposing new methods able to build a probability distribution on the space of sampled compatible models that represent preferences given by the DM.  In particular, we build some methodologies which assume that the DM takes into account a plurality of models specifications compatible with his preferences and each of them has a non zero contribution to the final alternatives' evaluation. In the framework of stochastic nature of the preference information, in this paper we aim to improve and overcome the previous results both in terms of reliability and efficiency. Our proposal takes into account a mixed approach which is based on the topological structure of the space of considered models (on which it is required the definition of a distance) together with a relatively simple inferential statistic treatment. In order to prove the goodness of our proposal we compared it to state of the art methods in this field performing a relatively large number of simulations.

The paper has the following structure. In Section \ref{prel} we present some preliminary notions and remarks useful for what follows together with the SMAA and the Subjective Stochastic Ordinal Regression (SSOR; \citealt{CorrenteEtAl2016KNOSYS}). Section \ref{indpro} describes the assumptions and the models underlying the proposed methods. In Section \ref{illex} we detail the computational experiment design and the statistical analysis used to test the efficiency of our proposal. Section \ref{Results} presents the obtained results. In Section \ref{SensitivitySection} we perform further computational experiments and a sensitivity analysis, while, in Section \ref{Comments} we discuss the results of all performed experiments. Finally, Section \ref{conc} contains conclusions and future directions of research.

%%%%%%%%%%%%%%%%%%%%%%%%%%%%%%%%%%%%%%%%%%%%%%%%%
\section{Notation and preliminaries}\label{prel}%
%%%%%%%%%%%%%%%%%%%%%%%%%%%%%%%%%%%%%%%%%%%%%%%%%
We consider a decision problem in which a set of alternatives $A=\left\{a,b,\ldots\right\}$\footnote{In the following, we shall denote the elements of $A$ by $a_j$ or $a$ interchangeably.}, is evaluated with respect to a set of criteria $G=\left\{g_1,\ldots,g_n\right\}$. Without loss of generality, we assume that, for all $g_i\in G$, ${g}_{i}:A\rightarrow {\rea}$ and it is of the gain type, that is, for all $a,b\in A$, ${g}_{i}(a)\geqslant {g}_{i}(b)$ means that $a$ is at least good as $b$ with respect to ${g}_{i}$. In what follows we assume as preference model a weighted sum
$$
U(a,\mathbf{w})=\displaystyle\sum_{i=1}^{n}w_i\cdot g_i(a)
$$
\noindent where $\mathbf{w}\in \mathbf{W}=\left\{(w_1,\ldots,w_n)\in[0,1]^{n}:\displaystyle\sum_{i=1}^{n} w_i=1\right\}$. In the following, to simplify the notation, we shall use $U(a)$ instead of $U(a,\mathbf{w})$ when it will not be necessary to refer to the weighted sum in particular. Pay attention to the fact that the procedure we are going to describe is not related to the choice of a weighted sum as preference model but it holds, more in general, considering a classical additive value function $U:A\rightarrow[0,1]$ assigning to each alternative $a\in A$ the value $\displaystyle\sum_{i=1}^{n}u_i(g_i(a))$, where $u_i$ is a marginal value function related to $g_i\in G$. Indeed, as will be clear later, it is enough considering a value function space ${\cal U}$ such that (i) it is possible to compute a distance between two elements in the space and, (ii) the preferences of the DM can be translated into linear constraints. 

Even if a weighted sum can be considered as the simplest MCDM model and, indeed, it is mainly used in practice just for this reason, its application implies the knowledge of the weights $w_i$ for all criteria $g_i$. In the following, we shall assume that the DM will use the indirect elicitation paradigm \citep{JacquetLagrezeSiskos2001} and, in particular, we distinguish between two cases:  % 

\begin{enumerate}
\item the DM does not provide certain preferences. In this case, it is common to consider all possible parameters specifications in the admissible space, that is, in our case, the simplex of all the weight vectors $\mathbf{W}$,
% in order to explore the preference relations determined by different preference parameters relations among alternatives (see for example \cite{lahdelma1998smaa});
\item the DM provides preference information in terms of a binary preference relation $\succsim$ on a subset of alternatives $A^R\subseteq A$ he knows relatively well, such that $a\succsim b$ iff $a$ is at least as good as $b$ and $a,b\in A^R$ (with $\succ$ and $\sim$ being the asymmetric and the symmetric part of relation $\succsim$, respectively), and another quaternary preference relation $\succsim^*$ on $A^R$ such that $(a,b)\succsim^{*}(c,d)$ iff the intensity of preference of $a$ over $b$ is at least as strong as the intensity of preference of $c$ over $d$, with $c\succ d$ and $a,b,c,d\in A^R.$ These preferences are then translated into constraints in terms of the assumed preference model. In particular, $a\succsim b$ iff $U(a)\geqslant U(b)$, while, $(a,b)\succsim^{*}(c,d)$ iff $U(a)-U(b)\geqslant U(c)-U(d)$ and $U(c)\geqslant U(d)+\varepsilon$ \citep{FigueiraGrecoSlowinski2009}, where $\varepsilon$ is an auxiliary variable, assumed greater than zero, used to convert strict inequalities into weak ones. To check if there exists at least one compatible model (in our case a weighted sum) one has to solve the following LP problem: 

%In this case, if the set of value functions compatible with the preference expressed by the DM is not empty, one can sample a well-distributed set of them or example, considering as in our case value functions expressed as weighted sum, one can use the Hit-And-Run (HAR) algorithm in the convex space of compatible value functions (Smith, 1984; van Valkenhoef et al. 2014) 
% a set of value functions compatible with such preferences can be computed starting from a feasible solution of the linear programming problem (LP) in eq.\ref{lpcert} and then applying an Hit and Run (HAR) algorithm (see \cite{smith1984efficient,van2014notes,kadzinski2016post}) to sample a set of a certain number of feasible solutions.\\
$$
\varepsilon^{*}=\max \varepsilon, \;\;\mbox{s.t.},\\
$$
\begin{equation}\label{lpcert}
\begin{array}{l}
\left.
\begin{array}{ll}
\;\;U(a)\geqslant U(b)+\varepsilon, \;\;\mbox{if $a\succ b$, for}\; a,b\in A^R\\[2mm]
\;\;U(a) = U(b), \;\;\mbox{if $a\sim b$, for}\; a,b\in A^R\\[2mm]
\left.
\begin{array}{l}
U(a)-U(b)\geqslant U(c)-U(d)+\varepsilon,\\[2mm]
U(c)-U(d)\geqslant \varepsilon,\\[2mm] 
\end{array}
\right\}\;\;\mbox{if $(a,b)\succ^*(c,d)$, for}\; a,b,c,d\in A^R\\[2mm]
\;\;U(a)=\displaystyle\sum_{i=1}^{n}w_i\cdot g_i(a), \;\;\forall a\in A,\\[2mm]
\;\;{w}_{i}\geqslant 0, \;\mbox{for all}\; i=1,\ldots,n,\\
\;\;\displaystyle\sum_{i=1}^{n}{w}_{i}=1.
\end{array}
\right\}E^{A^R}_{certain}
\end{array}
\end{equation}
If $E_{certain}^{A^R}$ is feasible and $\varepsilon^*>0$, then, the space of compatible value functions defined by constraints in $E_{certain}^{A^R}$, denoted by $\mathbf{W}_{certain}^{A^R}$, is not empty. In the opposite case ($E_{certain}^{A^R}$ is infeasible or $\varepsilon^*\leqslant 0$) there is not any compatible value function. Therefore, the cause of this inconsistency has to be checked and handled by using, for example, one of the two methods presented in \cite{MousseauEtAl2003} looking for several subsets of constraints of minimum cardinality that, once removed, restore the feasibility of $E^{A^{R}}_{certain}$. Presenting the DM with more than one of these sets is important since she/he can decide to remove the constraints on which she/he is less convinced.
\end{enumerate}
In both cases, since the spaces $\mathbf{W}$ and $\mathbf{W}_{certain}^{A^R}$ are composed of an infinite number of weight vectors and, consequently, an infinite number of corresponding weighted sum exists, one can give information regarding the alternatives at hand by considering a good approximation of these spaces. For such a reason, a sampling from these spaces has to be performed, and, in these cases, an algorithm such as Hit-And-Run (HAR; \citealt{Smith1984,TervonenEtAl2013}) can be used.

%%%%%%%%%%%%%%%%%%%%%%%%%%%%%%%%%%%%%%%%%%%%%%%%%%%%%%%%%%%%%%%%%%%%%%%%%%%%%%
\subsection{Stochastic Multicriteria Acceptability Analysis}\label{smaaintro}%
%%%%%%%%%%%%%%%%%%%%%%%%%%%%%%%%%%%%%%%%%%%%%%%%%%%%%%%%%%%%%%%%%%%%%%%%%%%%%%
Assuming a certain preference function (in our case the weighted sum  $U(a,\mathbf{w})=\displaystyle\sum_{i=1}^{n}{w}_{i}{g}_{i}(a)$), Stochastic Multicriteria Acceptability Analysis (SMAA) methods \citep{LahdelmaHokkanenSalminen1998,PelissariEtAl2020} provide robust information in statistical terms by considering as preference information a probability distribution function ${f}_{W}$ over $\mathbf{W}$ (in case 1. above) or a probability distribution function ${f}_{W_{certain}^{A^{R}}}$ over $\mathbf{W}_{certain}^{A^{R}}$ (in case 2. above) and a probability distribution function ${f}_{\chi}$ over the evaluation space $\mathbf{X}\subseteq {\rea}^{m\times n}$ being composed of all matrices $[{g}_{i}({a}_{j})]$ with $g_i\in G$ and $a_j\in A$. In our context, we shall assume that the evaluations will be deterministic and, therefore, we shall not consider ${f}_{\chi}$. Moreover, to simplify the notation, when it will not be necessary to distinguish between the two weight vector spaces $\mathbf{W}$ and $\mathbf{W}_{certain}^{A^{R}}$, we shall use $\mathbf{W}$.   

% once chosen a specific form for the value function, for example the linear value function $U(a,\mathbf{w})=\sum_{i=1}^{n}{w}_{i}{g}_{i}(a),$ where ${w}_{i}\geqslant 0$ for all $g_i\in G$ and $\sum_{i=1}^{n}{w}_{i}=1,$ if the evaluations are collected in a matrix $[{g}_{i}({a}_{j})]$ with $g_i\in G$ and $a_j\in A$, the preference information is composed of two probability distribution functions ${f}_{\chi}$ and ${f}_{W}$ defined over the evaluation space $\mathbf{X}\subseteq {\rea}^{m\times n}$ and over the weight space\\
%$$\displaystyle\mathbf{W}=\{({w}_{1},\ldots,{w}_{n}): {w}_{i}\geq 0\; \text{and}\;\displaystyle\sum_{i=1}^{n}{w}_{i}=1 \},$$ \\
%respectively, being $\chi:\mathbf{X}\rightarrow\rea$ and $W:\mathbf{W}\rightarrow\rea$ two measurable functions. 
%This paper is focused on the induction of a probability distribution function ${f}_{W}$ over the space of sampled compatible value functions approximating  the ``assumed'' probability distribution function of the DM denoted from now on by ${f}^{DM}_{}$.

Knowing the probability distribution ${f}_{W}$ over the weight space $\mathbf{W}$, the following indices are computed in SMAA-2 \citep{lahdelma2001smaa}:
\begin{itemize}
\item \textit{The rank acceptability index}:
\begin{equation}\label{RAIIntegral}
{b}^{r}(a)=\int_{\mathbf{w}\in {\mathbf{W}}_{}^{r}(a)}f_{W}(\mathbf{w})\;d\mathbf{w}
\end{equation}
\noindent where ${\mathbf{W}}_{}^{r}(a)=\left\{\mathbf{w}\in \mathbf{W}:rank(a,\mathbf{w})=r\right\},$ and $rank(a,\mathbf{w})=1+\displaystyle\sum_{x\in A\setminus \{a\}}\rho\left(U(x,\mathbf{w})>U(a,\mathbf{w})\right),$ with $\rho(false)=0$ and $\rho(true)=1$. The range of this index is $[0,1],$ meaning that the greater its value, the greater is the probability that the alternative $a$ achieves the rank $r$;

%Note that, since we consider deterministic evaluations ${g}_{i}(a)$ of considered alternatives $a \in A$ on criteria ${g}_{i} \in G$, we have that ${f}_{\chi}={I}_{\{\chi = \xi \}}$ and\\
%$$
%{b}^{r}_{a}=\int_{\mathbf{w}\in {\mathbf{W}}_{}^{r}(a,\xi)}f_{W}(\mathbf{w})\;d\mathbf{w}.
%$$

%Clearly, the best alternatives are those having rank acceptability index greater than zero for the first positions and rank acceptability index close to zero for the lower positions and it is within the range $[0,1]$.

\item \textit{The pairwise winning index} \citep{LeskinenEtAl2006}:
\begin{equation}\label{PWIIntegral}
p(a,b)=\int_{\mathbf{w}\in \mathbf{W}(a,b)} {f}_{W}(\mathbf{w}) \; d\mathbf{w}
\end{equation}
\noindent where $\mathbf{W}(a,b)=\{\mathbf{w}\in\mathbf{W}: U(a,\mathbf{w})> U(b,\mathbf{w})\}$. $p(a,b)$ represents the probability with which an alternative $a$ is strictly preferred to an alternative $b$ in $\mathbf{W}.$ This index is also within the range $[0,1]$ and the greater $p(a,b)$, the more $a$ is preferred to $b$. 
\end{itemize}

Usually, the rank acceptability indices and the pairwise winning indices are collected into matrices $RAI$ and $PWI$ ($RAI=[{b}^{r}(a)],$ $PWI=[{p}(a,b)]$) for a global overview. Clearly, their values depend on $f_{W}$. Consequently, we shall use the symbols $RAI_{f_{W}}$ and $PWI_{f_{W}}$ to denote the $RAI$ and $PWI$ matrices obtained assuming the probability distribution $f_{W}$ over $\mathbf{W}$. For many convenient reasons, $f_W$ is often chosen among some parametric distribution functions. In many works (see, for example, \citealt{TervonenFigueira2008,PelissariEtAl2020}), $f_W$ follows a discrete uniform distribution when no preference information is available. Otherwise, $f_W \sim N (\mu,\Sigma)$ (truncated) when some preference information suggests it (see, for example, \citealt{LahdelmaMakkonenSalminen2006,KangasEtAl2006,LahdelmaSalminen2010}). 

Differently from the works mentioned above, in this paper, we do not assume the knowledge of the probability distribution $f_{W}$ but we want to propose some methods aiming to infer this probability on the basis of the preference information provided by the DM. For such a reason, we shall denote this probability distribution by $f_{DM}$.  

%Observe that since we consider deterministic evaluations $g_i(a), a \in A g_i \in G$, the pairwise winning index can be reformulated as 
%$$P_{a,b}=/int_{mathbf{w} \in W: U(a,\mathbf{w})>U(b,\mathbf{w})} f_W(\mathbf{w})d \mathbf{w}$.

Let us conclude this section by observing that in SMAA we do not consider the whole infinite space $\mathbf{W}$, but a well-distributed sample of weight vectors in it. We shall denote the sample of weight vectors in $\mathbf{W}$ by $\mathbf{\Omega}$ and, consequently, $\mathbf{\Omega}\subseteq\mathbf{W}$. In this context, each weight vector $\mathbf{w}\in\mathbf{\Omega}$ has a mass $p(\mathbf{w})\in\left[0,1\right]$ such that $\displaystyle\sum_{\mathbf{w}\in \mathbf{\Omega}}p(\mathbf{w})=1$. Taking into account the space $\mathbf{\Omega}$ of sampled weight vectors, for all $a,b \in A$, the rank acceptability index $b^r(a)$ and the pairwise winning index $p(a,b)$ in eqs. (\ref{RAIIntegral}) and (\ref{PWIIntegral}), respectively, can be formulated as follows:
\[
b^r(a)=\displaystyle\sum_{\mathbf{w}\in\mathbf{\Omega}_r(a)}p(\mathbf{w})\qquad\mbox{and}\qquad p(a,b)=\displaystyle\sum_{\mathbf{w}\in\mathbf{\Omega}(a,b)}p(\mathbf{w})
\]
\noindent where $\mathbf{\Omega}_r(a)=\left\{\mathbf{w}\in\mathbf{\Omega}:\;rank(a,\mathbf{w})=r\right\}$ and $\mathbf{\Omega}(a,b)=\{\mathbf{w} \in\mathbf{\Omega}: U(a,\mathbf{w}) >U(b,\mathbf{w})\}$.

For this reason, we can reformulate the aim of this paper saying that we want to introduce some methods to build a probability distribution on $\mathbf{\Omega}$ on the basis of certain and uncertain preferences given by the DM. 

% From this point of view, the aim of this paper is to propose some methods to estimate ${f}^{DM}_{W}$.

%%%%%%%%%%%%%%%%%%%%%%%%%%%%%%%%%%%%%%%%%%%%%%%%%%%%%%%%%%%%%%%%%%%%%%%
\subsection{Subjective Stochastic Ordinal Regression}\label{ssorlabel}%
%%%%%%%%%%%%%%%%%%%%%%%%%%%%%%%%%%%%%%%%%%%%%%%%%%%%%%%%%%%%%%%%%%%%%%%
The Subjective Stochastic Ordinal Regression (see \citealt{CorrenteEtAl2016KNOSYS}) aims to define a probability distribution $\mathbf{p}_{DM}$ over $\mathbf{\Omega}$ assigning, therefore, a mass $p(U_t)$ to each value function $U_t\in\mathbf{\Omega}$. For this reason, it is assumed that a DM provides uncertain preferences in terms of pairwise comparisons of alternatives such as ``\textit{the preference of $a$ over $b$ is at least as credible as the preference of $b$ over $a$}'' (denoted by $a{\succsim}_{Pr}b$) or in terms of intensity of preference such as ``\textit{the preference of $a$ over $b$ is at least as credible as the preference of $c$ over $d$}'' ($(a,b){\succsim}^{*}_{Pr}(c,d)$). Note that, from a probabilistic point of view, the uncertain information $a{\succsim}_{Pr}b,$ can be read as ``the probability that $a$ is at least as good as $b$ is not lower than the probability that $b$ is at least as good as $a,$'' while, the uncertain information $(a,b){\succsim}^{*}_{Pr}(c,d)$ can be read as ``the probability that $a$ is at least as good as $b$ is not lower than the probability that $c$ is at least as good as $d$''. These preference relations are used as constraints to induce the DM subjective probability measure on the set of sampled compatible value functions $\mathbf{\Omega}$.

From a computational point of view, following \cite{CorrenteEtAl2016KNOSYS}, to find a probability distribution compatible with the uncertain preference information provided by the DM one has to solve the following LP problem: % , and, if it is feasible for an $\varepsilon>0,$ to compute the barycenter of a sampling of measures ($\Omega$) compatible with the constraints in ${E}^{SSOR}.$
$$
\varepsilon_{L}^{*}=\max\varepsilon, \;\;\mbox{s.t.},\\
$$
\begin{equation}\label{ssorlpp}
\begin{array}{l}
%\;\;\;\max\varepsilon, \;\;\mbox{subject to}\\[3mm]
\left.
\begin{array}{l}
\displaystyle\sum_{U_t\in\mathbf{\Omega}:\; U_t(a)>U_t(b)}p(U_t)\geqslant \sum_{U_t\in\mathbf{\Omega}:\; U_t(b)>U_t(a)}p(U_t)+\varepsilon, \;\;\;\mbox{if $a\;\succ_{P_r}\;b$},\\[4mm]
\displaystyle\sum_{U_t\in\mathbf{\Omega}:\; U_t(a)>U_t(b)}p(U_t)\geqslant \sum_{U_t\in\mathbf{\Omega}:\; U_t(c)>U_t(d)}p(U_t)+\varepsilon, \;\;\;\mbox{if $(a,b)\;{\succ}^{*}_{Pr}\;(c,d)$},\\[4mm]
\displaystyle\sum_{t=1}^{|\mathbf{\Omega}|}p(U_t)=1,\\[5mm]
p(U_t)\geqslant  0, t=1,\ldots,|\mathbf{\Omega}|
\end{array}
\right\}{E}^{SSOR}
\end{array}
\end{equation}

\noindent where $\succ_{P_r}$ and $\succ^{*}_{P_r}$ denote the asymmetric parts of $\succsim_{P_{r}}$ and $\succ_{P_r}$, respectively, and $p(U_t)$, with $U_t\in\mathbf{\Omega}$, is the ``mass'' attached to the sampled value function $U_t$. As already underlined in the previous section, $\mathbf{\Omega}$ represents a sample of value functions compatible with the preferences provided by the DM (if any) and these functions do not have to be necessarily weighted sum. If ${E}^{SSOR}$ is feasible and $\varepsilon_{L}^{*}>0$, there is at least one probability distribution over $\mathbf{\Omega}$ compatible with the uncertain preferences given by the DM. Since, in this case, more than one probability distribution compatible with the preferences given by the DM exists, \cite{CorrenteEtAl2016KNOSYS} define a probabilistic necessary and a probabilistic possible preference relation taking into account the whole space of probability distribution defined by constraints in $E^{SSOR}$ (the interested reader can find more details in \citealt{CorrenteEtAl2016KNOSYS}). However, they suggest also giving information in statistical terms by sampling a certain number of probability distributions from the space defined by constraints in $E^{SSOR}$ and, then, considering their barycenter $\mathbf{p}^{*}$ as a representative of them. On the basis of this barycenter, a binary relation $\succsim^{R}_{L}$ was defined so that $a\succsim^{R}_{L} b$ iff, considering $\mathbf{p}^{*}$, the probability that $a$ is preferred to $b$ ($PWI_{\mathbf{p}^*}(a,b)$) is not lower than the probability that $b$ is preferred to $a$ ($PWI_{\mathbf{p}^*}(b,a)$).

Let us underline that the probability distribution obtained by the $SSOR$ assigns the masses neglecting any regularity requirement on the probability distribution on the set of sampled value functions $\mathbf{\Omega}$ so that this methodology can attach very different masses to value functions relatively close in $\mathbf{\Omega}$. In the next section, we shall present some methods to build a probability distribution on $\mathbf{\Omega}$ that use the same information as SSOR but that assign the probability mass to each sampled value function depending on their distance from a reference value function $U_{ref}$. 

%%%%%%%%%%%%%%%%%%%%%%%%%%%%%%%%%%%%%%%%%%%%%%%%%%%%%%%%%%%%%%%%%%%%%%%%%%%%%%%%%%%%%%%%%%%%%%%%%%%%%
\section{Inferring a probability density function over the space of compatible models}\label{indpro}%
%%%%%%%%%%%%%%%%%%%%%%%%%%%%%%%%%%%%%%%%%%%%%%%%%%%%%%%%%%%%%%%%%%%%%%%%%%%%%%%%%%%%%%%%%%%%%%%%%%%%%
In this section, we present a methodology to infer a probability distribution function $\mathbf{p}=\left[p(U_t):\;U_t\in\mathbf{\Omega}\right]$ over a sampling of compatible models $\mathbf{\Omega}$ (in our case a probability distribution $\mathbf{p}$ on the sample $\mathbf{\Omega} \subseteq \mathbf{W}$ composed of weight vectors). Note that, in general terms, the methods we propose can be applied in case the assumed preference model is defined on a space ${\cal U}$ equipped with a distance and such that, rewriting conveniently the last three constraints in (\ref{lpcert}), the programming problem remains linear. For example:
\begin{itemize}
\item if the assumed preference model is an additive value function of the type $U(a)=\displaystyle\sum_{i=1}^{n}u_i\left(g_{i}(a)\right)$, the last three constraints in (\ref{lpcert}) should be replaced by the following 
$$
\left.
\begin{array}{l}
U(a)=\displaystyle\sum_{i=1}^{n}u_i\left(g_{i}(a)\right),\;\mbox{for all}\;a\in A,\\[2mm]
\displaystyle\sum_{i=1}^{n}u_i\left(x_i^{*}\right)=1,\\[2mm]
u_i\left(g_i(a)\right)\geqslant u_i\left(g_i(b)\right),\;\mbox{iff}\; g_i(a)\geqslant g_i(b), \;\mbox{for all}\; i=1,\ldots,n,\\[2mm]
u_i\left(x_{i,*}\right)=0\;\mbox{for all $i=1,\ldots,n$},\\[2mm]
\end{array}
\right\}
$$
\noindent where, for all $i=1,\ldots,n$, $x_{i,*}=\displaystyle\min_{a\in A}g_i(a)$ and $x_{i}^{*}=\displaystyle\max_{a\in A}g_i(a)$, respectively,
%$a_*,a^*\in A$ are such that, for all $i=1,\ldots,n,$ $g_i(a_*)=\displaystyle\min_{a\in A}g_i(a)$ and $g_i(a^*)=\displaystyle\max_{a\in A}g_i(a)$, respectively, 
\item if the preference models is the 2-additive Choquet integral \citep{Choquet1953,Grabisch1996,Grabisch1997} $U(a)=Ch_{\mu}(a)=\displaystyle\sum_{g_i\in G}g_i(a)m\left(\{g_i\}\right)+\sum_{\{g_{i_1},g_{i_2}\}\subseteq G}\min\{g_{i_1}(a),g_{i_2}(a)\}m\left(\{g_{i_1},g_{i_2}\}\right)$, then, the last three constraints in (\ref{lpcert}) should be replaced by the following ones
$$
\left.
\begin{array}{l}
\;\;U(a)=Ch_{\mu}(a)=\displaystyle\sum_{g_i\in G}g_i(a)m\left(\{g_i\}\right)+\sum_{\{g_{i_1},g_{i_2}\}\subseteq G}\min\{g_{i_1}(a),g_{i_2}(a)\}m\left(\{g_{i_1},g_{i_2}\}\right),\;\mbox{for all}\;a\in A,\\[2mm]
\left.
\begin{array}{l}
m\left(\{g_i\}\right)\geqslant 0,\\[2mm]
m\left(\{g_i\}\right)+\displaystyle\sum_{g_{i_{1}}\in T}m\left(\{g_{i},g_{i_1}\}\right)\geqslant 0,
\end{array}
\right\}\mbox{for all}\; g_i\in G\; \mbox{and for all}\; T\subseteq G\setminus\{g_i\}\\[2mm]
\displaystyle\sum_{g_i\in G}m\left(\{g_i\}\right)+\sum_{\{g_{i_1},g_{i_2}\}\subseteq G}m\left(\{g_{i_1},g_{i_2}\}\right)=1
\end{array}
\right\}
$$  
\noindent where $m:2^{G}\rightarrow\rea$ is a set function such that $m(T)=0$ for all $T\subseteq G$ with $|T|>2$. In fact, the model we are proposing can be applied also when constraints in (\ref{lpcert}) are non-linear, for example, when the preference of the DM can be represented by a multiplicative utility function \citep{Keeney1974, KeeneyRaiffa1976}. In these cases, one can apply an acceptance/rejection method (see e.g. Section 2.5.1 in \citealt{RubinsteinKroese2016}), based on Hit-And-Run, sampling in a convex set containing the set of parameters of the non-linear value function considered.
\end{itemize}

%%%%%%%%%%%%%%%%%%%%%%%%%%%%%%%%%%%%%%%%%
\subsection{Basic proposal\label{mainp}}%
%%%%%%%%%%%%%%%%%%%%%%%%%%%%%%%%%%%%%%%%%
Starting from the work presented in \cite{CorrenteEtAl2016KNOSYS}, we propose to define a probability distribution function over $\mathbf{\Omega}$ with a tendency preference (a mode) which we refer to as a reference model ${U}_{ref}\in{\cal U}$ where ${\cal U}$ denotes the set of models compatible with the preferences given by the DM (indeed, assuming that the preference model is a weighted sum, ${\cal U}=\{U(a,\mathbf{w}): \mathbf{w} \in \mathbf{W}\}$ if the DM did not provide any preference information, while ${\cal U}=\{U(a,\mathbf{w}): \mathbf{w} \in \mathbf{W}^{A^R}_{certain}\}$ if the DM expressed some preferences). Such choice has two motivations: (i) the direction given by the uncertain information (see Section \ref{ssorlabel}) suggests that there are some models which better represent the preferences provided by the DM, so, it can be reasonable to take a single ``average model'' representing the typical attitude of the DM; (ii) the induction of a probability distribution with a mode in the space of compatible models allows to easily ``move'' the main body of the distribution, i.e. the largest probability, in the direction indicated by the preferences expressed by the DM. With this in mind, we propose to compute the reference model ${U}_{ref}$ as the barycenter of the set of sampled value functions $\mathbf{\Omega}$. \\
The barycenter is computed as follows: 
\begin{description}
\item[Step 1)] Solve the LP problem (\ref{lpcert}). If $E_{certain}^{A^R}$ is feasible and $\varepsilon^*>0$, go to step 2) below. If this is not the case ($E^{A^R}_{certain}$ is not feasible or $\varepsilon^*\leqslant 0$), we have to interact with the DM in order to remove a minimal number of constraints related to preference information and make  $E^{A^R}_{certain}$ feasible, using one of the two procedures proposed in \cite{MousseauEtAl2003} and recalled in Section \ref{prel},
\item[Step 2)] Sample a certain number of compatible value functions using the HAR algorithm defining the set $\mathbf{\Omega}$,
\item[Step 3)] Since, in our case, each compatible value function $U_t\in \mathbf{\Omega}$ is the weighted sum related to the weight vector $\mathbf{w}_t=\left(w_1^{t},\ldots,w_n^{t}\right)$, the barycenter of $\mathbf{\Omega}$ is also a weight vector $\mathbf{w}_{Bar}=\left(w_1^{Bar},\ldots,w_n^{Bar}\right)$ where $w_{i}^{Bar}=\frac{1}{|\mathbf{\Omega}|}\displaystyle\sum_{t=1}^{|\mathbf{\Omega}|}w_i^{t}$, for all $i=1,\ldots,n$. 
\end{description}
Following what has been said at the beginning of this section, if the assumed preference model is an additive value function or a 2-additive Choquet integral, then, the barycenter can be computed as explained in Steps 1)-3) above. The only slight difference is related to Step 3): \\
\begin{itemize}
\item if $U_t\in\mathbf{\Omega}$ is an additive value function, $U_t=\left[u^{t}_i(g_i(a))\right]_{\substack{i=1,\ldots,n \\ a\in A}}$ and, consequently, the barycenter is an additive value function $U_{Bar}=\left[u^{Bar}_i(g_i(a))\right]_{\substack{i=1,\ldots,n \\ a\in A}}$ such that $u_i^{Bar}(g_i(a))=\frac{1}{|\mathbf{\Omega}|}\displaystyle\sum_{t=1}^{|\mathbf{\Omega}|}u^{t}_i(g_i(a)),$
\item if $U_t\in\mathbf{\Omega}$ is a 2-additive value function, $U_t=\left[\left[m^{t}(\{g_{i}\})\right]_{g_i\in G},\left[m^{t}(\{g_{i_1},g_{i_2}\})\right]_{\{g_{i_1},g_{i_2}\}\subseteq G}\right]$ and, consequently, the barycenter is a 2-additive Choquet integral $$
U_{Bar}=\left[\left[m^{Bar}(\{g_{i}\})\right]_{i=1,\ldots,n},\left[m^{Bar}(\{g_{i_1},g_{i_2}\})\right]_{\{g_{i_1},g_{i_2}\}\subseteq G}\right]
$$ 
\noindent such that $\displaystyle m^{Bar}\left(\{g_i\}\right)=\frac{1}{|\mathbf{\Omega}|}\sum_{t=1}^{|\mathbf{\Omega}|}m^{t}\left(\{g_i\}\right)$ and $m^{Bar}\left(\{g_{i_1},g_{i_2}\}\right)=\displaystyle\frac{1}{|\mathbf{\Omega}|}\sum_{t=1}^{|\mathbf{\Omega}|}m^{t}\left(\{g_{i_1},g_{i_2}\}\right)$ for all $g_{i}\in G$ and for all $\left\{g_{i_1},g_{i_2}\right\}\subseteq G$.
\end{itemize}

We search for a probability distribution function $\mathbf{p}$ over $\mathbf{\Omega}$ which assigns a mass $p(U_t)$ to each preference model $U_t\in\mathbf{\Omega}$ that, differently from SSOR, is decreasing with respect to the distance of $U_t$ from $U_{ref}$. This procedure requires therefore the introduction of a distance in the space $\mathbf{\Omega}$. If the preference model can be represented by means of a real valued vector of preferential parameters, the Euclidean norm and its induced metric (as well as each equivalent norm) can be considered\footnote{For example, if the preference model is a weighted sum, then, this function can be represented by the vector $\mathbf{w}=\left[w_1,\ldots,w_n\right]$; if the preference function is a general non-monotonic additive value function, then, it can be represented by the vector $\left[u_i(g_i(a))\right]_{\substack{g_i\in G\\ a\in A}}$, and so on.}. 

All previous considerations lead to the search for a probability distribution function $\mathbf{p}$ over $\mathbf{\Omega}$ in which, given $U_k,U_h\in\mathbf{\Omega}$, ${p}({U}_{k})\geqslant {p}({U}_{h})$ iff $d({U}_{k},{U}_{ref})\leqslant d({U}_{h},{U}_{ref}),$ where $d\left({U}_{h},{U}_{ref}\right)$ is the distance between the vectors representing the value functions $U_h$ and $U_{ref}$. For example, if the value function is a weighted sum and, therefore, $U_h$ is represented by the vector $\mathbf{w}_{h}=\left[w_{h,1},\ldots,w_{h,n}\right]$, then, one can have

$$d({U}_{h},{U}_{ref})=\sqrt{\displaystyle\sum_{i=1}^{n}\left(w_{h,i}-w_{ref,i}\right)^2}.$$

Assuming the existence of some uncertain preference information as those considered in \cite{CorrenteEtAl2016KNOSYS}, to check if there exists at least one probability distribution over $\mathbf{\Omega}$, one has to solve the following LP problem (\ref{lpp1}):
$$
\varepsilon_{ACG}^{*}=\max{\varepsilon},\; \mbox{s.t.}\\
$$
\begin{equation}\label{lpp1}
\left.
\begin{array}{l}
\displaystyle\sum_{U_t\in\mathbf{\Omega}:\;U_t(a)\geqslant U_t(b)} {p}({U}_{t})\geqslant \sum_{U_t\in\mathbf{\Omega}:\;U_t(a)\leqslant U_t(b)} {p}({U}_{t}) + \varepsilon, \;\;\mbox{if $a{\succ}_{Pr} b$},\\[6mm]
\displaystyle\sum_{U_t\in\mathbf{\Omega}:\;U_t(a)\geqslant U_t(b)} {p}({U}_{t})\geqslant \sum_{U_t\in\mathbf{\Omega}:\;U_t(c)\geqslant U_t(d)} {p}({U}_{t}) + \varepsilon, \;\;\mbox{if $(a,b){\succ}^{*}_{Pr}(c,d)$},\\[6mm]
{p}({U}_{k})\geqslant {p}({U}_{h}), \mbox{ iff $d({U}_{k},{U}_{ref})\leqslant d({U}_{h},{U}_{ref})$},\; \mbox{for all ${U}_{k},{U}_{h}\in\mathbf{\Omega}, $} \\[2mm]
\displaystyle\sum^{|\mathbf{\Omega}|}_{t=1} {p}({U}_{t})= 1,\\[6mm]
p({U}_{t})\geqslant 0, \mbox{ for all ${U}_{t}\in \mathbf{\Omega}$}. \\[2mm]
% \end{array}
% \right \}
\end{array}
\right \}{E}^{ACG}
\end{equation}

\noindent If $E^{ACG}$ is feasible and $\varepsilon_{ACG}^{*}>0$ then, there exists at least one probability distribution having the previous characteristics and, therefore, we can sample by the HAR algorithm a certain number of probability distributions from the space defined by constraints in $E^{ACG}$ computing then its barycenter. The barycenter is therefore used to compute the rank acceptability index of each alternative for each rank position as well as the pairwise winning index between each ordered pair of alternatives.  If $E^{ACG}$ is not feasible or $\varepsilon_{ACG}^{*}\leqslant 0$, one can check for the causes as discussed in Section \ref{prel}. 

\subsection{Estimation of a parsimonious model}\label{piecep}%
%%%%%%%%%%%%%%%%%%%%%%%%%%%%%%%%%%%%%%%%%%%%%%%%%%%%%%%%%%%%%%
As explained in the previous section, in the $SSOR^{ACG}$ method the number of parameters to be estimated increases with the size of the sample since one probability mass is assigned to each sampled compatible model. For this reason, in the following, we propose a simple and parsimonious procedure to assign a probability mass to each sampled model assuming that the probability function is a piecewise linear value function defined by means of a certain number of reference points the masses of which will therefore be the only unknowns of our problem. The probability mass of the sampled value functions which are not given as reference points will be obtained by linear interpolation as explained in the following. Of course, also in this case, we assume that the probability mass attached to each model $U_h$ will be a non-increasing function of its distance from a reference model $d(U_h,U_{ref})=d^{ref}\left(U_h\right)$, that is, $p(U_h)=p\left(d^{ref}\left(U_h\right)\right)$.

%Given the previous considerations and assumptions on the shape of the distribution function, it is reasonable to assume that the probability assignment can be done according with a continuous weighting function $p$ which is decreasing with respect to the distance from the reference model. However, from a computational point of view, the 
%% ${U}_{q}=\displaystyle\operatorname*{arg\,max}_{{U}_{}\in\Omega}{d}^{r}({U}_{})$
Let us fix $q$ reference distances $d_1,...d_q$  into the interval $\mathcal{I}({U}_{ref},\mathbf{\Omega})=\left[0,\displaystyle\max_{h=1,\ldots,|\mathbf{\Omega}|}d^{ref}\left(U_h\right)\right]$, such that $d_{k-1}<d_k$ for all $k=2,\ldots,q $, $d_1=0$ and $\displaystyle d_q=\max_{h\in \{1,\ldots,|\mathbf{\Omega}|\}}d^{ref}\left(U_h\right)$. We aim to infer a probability mass function $p:\mathcal{I}({U}_{ref},\mathbf{\Omega})\rightarrow [0,1]$ such that 
if ${d}^{ref}({U}_{t})\in [{d}_{k-1},{d}_{k}],$ with $k=2,\dots,q$, we have that\\ 
$$
p\left(U_t\right)=p({d}^{ref}({U}_{t})) = p({d}_{k-1}) \frac{{d}_{k}-{d}^{ref}({U}_{t})}{{d}_{k} - {d}_{k-1}} + p({d}_{k})\frac{{d}^{ref}({U}_{t})-{d}_{k-1}}{{d}_{k} - {d}_{k-1}}.
$$
Knowing $p(d_k)$ for all $k=1,\ldots,q$, one can then compute the probability mass attached to all functions in $\mathbf{\Omega}$ depending on their distance from $U_{ref}$ and, consequently, from the interval to which this distance belongs. Note that, the assumption on the probability mass being a decreasing function of the distance from $U_{ref}$ requires that $p({d}^{ref}({U}_{k}))\geqslant p({d}^{ref}({U}_{h}))$ for each $U_h,U_k\in\mathbf{\Omega}$ such that ${d}^{ref}({U}_{k})\leqslant {d}^{ref}({U}_{h}),$ which is guaranteed by imposing that  $p({d}_{k-1})\geqslant p({d}_{k})$ for all $h=2,\ldots,q$. To check if there exists at least one piecewise linear function defined by the vector $[p(d_1),\ldots,p(d_q) ]$ compatible with the provided preference, one has to solve the following LP problem (\ref{lpp2}):\\
$$
\varepsilon_{pl}^{*}=\max{\varepsilon},\; \mbox{s.t.}\\
$$
\begin{equation}\label{lpp2}
\left.
\begin{array}{l}
p\left(U_t\right)= p({d}_{k-1}) \frac{{d}_{k}-{d}^{ref}({U}_{t})}{{d}_{k} - {d}_{k-1}} + p({d}_{k})\frac{{d}^{ref}({U}_{t})-{d}_{k-1}}{{d}_{k} - {d}_{k-1}}, \;\;\mbox{if}\;\;{d}^{ref}({U}_{t})\in [{d}_{k-1},{d}_{k}], \\[4mm]
\displaystyle\sum_{U_t\in\mathbf{\Omega}:\;U_t(a)\geqslant U_t(b)}p({U}_{t})\geqslant \displaystyle\sum_{t:\;U_t(a)\leqslant U_t(b)}p({U}_{t}) + \varepsilon, \;\;\mbox{if $a{\succ}_{Pr} b$},\\[6mm]
\displaystyle\sum_{U_t\in\mathbf{\Omega}:\;U_t(a)\geqslant U_t(b)}p({U}_{t})\geqslant \sum_{U_t\in\mathbf{\Omega}:\;U_t(c)\geqslant U_t(d)}p({U}_{t}) + \varepsilon,\;\;\mbox{if $(a,b){\succ}^{*}_{Pr}(c,d)$},\\[6mm]
p({d}_{k-1})\geqslant p({d}_{k})\geqslant 0, \; \mbox{for all}\; k=2,\ldots,q, \\[2mm]
%\left.
%\begin{array}{l}
\displaystyle\sum^{|\mathbf{\Omega}|}_{t=1}p({U}_{t}) = 1.  \\[4mm]
%\end{array}
%\right\}
\end{array}
\right\}{E}_{pl}^{ACG}
\end{equation}

If $E_{pl}^{ACG}$ is feasible and $\varepsilon^{*}_{pl}>0$ then, there is at least one piecewise linear probability distribution on the set of sampled value functions compatible with the uncertain preference given by the DM. Sampling a certain number of probability distributions in the space defined by the constraints in $E^{ACG}_{pl}$ by using the HAR algorithm we can compute its barycenter. The so-obtained barycenter can then be used to compute the RAIs and PWIs. We refer to this method as $SSOR_{pl}^{ACG}$.

\subsection{Parametric estimation}\label{paramp}%
%%%%%%%%%%%%%%%%%%%%%%%%%%%%%%%%%%%%%%%%%%%%%%%%%
The proposals presented in Sections \ref{mainp} and \ref{piecep} are based on a non-parametric approach in the sense that we do not assume ``ex-ante'' a specific shape of the probability distribution function over $\mathbf{\Omega}.$ In this section, instead, we aim to fit the distribution function defined over the set of all sampled compatible models by using some specific parametric probability function. As in the previous proposals, we assume that such functions are expressed in terms of distance from a reference model ${U}_{ref}$ and, in particular, the probability mass attached to each compatible model $U_h\in\mathbf{\Omega}$ is decreasing with its distance from $U_{ref}$, that is, $d^{ref}\left(U_h\right)=d\left(U_{ref},U_h\right)$. 

We start by considering the following two alternative probability functions for a  generic model ${U}_{h}\in\mathbf{\Omega}$:\\
\begin{equation}\label{nlnor}
{p}_{Nor}({U}_{h}) = \frac{{e}^{-\frac{{{d}^{ref}\left({U}_{h}\right)}^{2}}{2{\lambda}^{2}}}}{\displaystyle\sum_{h=1}^{|\mathbf{\Omega}|}{e}^{-\frac{{{d}^{ref}\left({U}_{h}\right)}^{2}}{2{\lambda}^{2}}}},
\end{equation}
%\begin{equation}\label{nlnor}
%{p}_{Nor}({U}_{i}) = \frac{\frac{\sqrt{2}}{\lambda\sqrt{\pi}}{e}^{-\frac{{{d}^{r}({U}_{i})}^{2}}{2{\lambda}^{2}}}}{\displaystyle\sum_{j=1}^{|\Omega|}\frac{\sqrt{2}}{\lambda\sqrt{\pi}}{e}^{-\frac{{{d}^{r}({U}_{i})}^{2}}{2{\lambda}^{2}}}} = \frac{{e}^{-\frac{{{d}^{r}({U}_{i})}^{2}}{2{\lambda}^{2}}}}{\displaystyle\sum_{j=1}^{|\Omega|}{e}^{-\frac{{{d}^{r}({U}_{i})}^{2}}{2{\lambda}^{2}}}},
%\end{equation}

\noindent and\\

\begin{equation}\label{nlexp}
{p}_{Exp}({U}_{h}) = \frac{{e}^{-\lambda {d}^{ref}\left({U}_{h}\right)}}{\displaystyle\sum_{h=1}^{|\mathbf{\Omega}|}{e}^{-\lambda {d}^{ref}\left({U}_{h}\right)}},
\end{equation}
%\begin{equation}\label{nlexp}
%{p}_{Exp}({U}_{i}) = \frac{\lambda {e}^{-\lambda {d}^{r}({U}_{i})}}{\displaystyle\sum_{j=1}^{|\Omega|}\lambda {e}^{-\lambda {d}^{r}({U}_{i})}} = \frac{{e}^{-\lambda {d}^{r}({U}_{i})}}{\displaystyle\sum_{j=1}^{|\Omega|}{e}^{-\lambda {d}^{r}({U}_{i})}}.
%\end{equation}
\noindent with $\lambda \geqslant 0$.

Now, we consider the mathematical programming problems (\ref{ppnor}) and (\ref{ppexp})

\begin{equation}\label{ppnor}
% \displaystyle\max{\sum_{U\in B}^{}x{e}^{-x{d}^{r}({U}_{})} - \sum_{U\in {B}^{c}}^{}x{e}^{-x{d}^{r}({U}_{})}},\; \mbox{s.t.}\\
\max_{}\displaystyle\sum_{{U}_{h}\in B}^{}{p}_{Nor}({U}_{h}),\; \mbox{s.t.}
\end{equation}
\begin{equation*}
\left.
%\begin{array}{l}
%\displaystyle af(x,\Omega)=1,\\% \sum_ {i=1}^{|\Omega|}x{e}^{-x{d}^{r}({U}_{i})}
\lambda\geqslant 0
%\end{array}
\right.% \}{E}^{GCAnl}
\end{equation*}

\begin{equation}\label{ppexp}
% \displaystyle\max{\sum_{U\in B}^{}x{e}^{-x{d}^{r}({U}_{})} - \sum_{U\in {B}^{c}}^{}x{e}^{-x{d}^{r}({U}_{})}},\; \mbox{s.t.}\\
\max_{}\displaystyle\sum_{{U}_{h}\in B}^{}{p}_{Exp}({U}_{h}),\; \mbox{s.t.}
\end{equation}
\begin{equation*}
\left.
%\begin{array}{l}
%\displaystyle af(x,\Omega)=1,\\% \sum_ {i=1}^{|\Omega|}x{e}^{-x{d}^{r}({U}_{i})}
\lambda\geqslant 0
%\end{array}
\right.% \}{E}^{GCAnl}
\end{equation*}

\noindent where $B=\{U_h\in \mathbf{\Omega}:\;{U}_{h}(a)\geqslant {U}_{h}(b),\;\forall\; (a,b)\in A\times A:\; a\succsim_{P_r}b\}$. Observe that these mathematical programming problems are well defined if $B$ is non-empty. If this is not the case, one can reformulate (\ref{ppnor}) and (\ref{ppexp}) replacing $B$ with $B^k=\{U_h \in \mathbf{\Omega}:\;|(a,b) \in A \times A\;:\; U_h(a)\;\geqslant\; U_h(b)$ and $a \succsim_{pr}b| \geqslant k\}$ for some value of $k$ for which $B^k$ is not empty. Let us observe that $k$ can be considered as the minimum number of pairwise comparisons on which the DM expressed his preferences that need to be represented by at least one compatible model. Finally, we compare the optimal objective values of the mathematical programming problems (\ref{ppnor}) and (\ref{ppexp}) and we consider the model which gives the greatest of the two. Note that, the proposal at hand differs from the previous not only from a parametric point of view but also from the fact that we use the probability distribution obtained as solution of the considered problems and not the barycenter of a sample of compatible distributions. %An example of this method, denoted as ${SSOR}^{GCAnl},$ is given in Example \ref{ex3}. 
\section{Validation proposal}\label{illex}%
In this section we present the results of some experiments performed to show the reliability of the methodologies aiming to build a probability distribution on the space of models compatible with the preferences given by the DM and described in Sections \ref{mainp}-\ref{paramp}. To compare the considered methods we mainly look at the difference between the $RAI$ and $PWI$ matrices generated by the probability distributions $p(\cdot)$ obtained with the different methods we are proposing. 

Assuming the existence of an artificial DM that ranks some alternatives evaluated on a few criteria, we aim to check which method, among the three we have proposed, is more able to restore the artificial DM's preferences. To simulate the artificial DM's preferences, we consider two different cases: 
\begin{itemize}
\item the artificial DM evaluates the alternatives using a single model $U^{DM}$ on which no assumption is done. It is only necessary that it is defined over a space ${\cal U}$ so that the distance between two value functions in this space can be computed. In this case, it is possible to define a binary relation $\succsim^{DM}$ over $A$ such that, for all $a,b \in A$, $a \succsim^{DM} b$ iff $U^{DM}(a) \geqslant U^{DM}(b)$;
\item the preferences of the DM are determined by a set of models $U_t\in\mathbf{\Omega}$, where $\mathbf{\Omega}$ is a set of well-distributed models in ${\cal U}$ on which a probability distribution $\mathbf{p}^{DM}=\left[p^{DM}(U_t): U_t\in\mathbf{\Omega}\right]$ is considered. In this case we can define the binary relation $\succsim^{DM}$ over $A$ so that for all $a,b \in A,\;a \succsim^{DM} b$ iff 
$$
\sum_{U_h \in\mathbf{\Omega}:\; U_h(a) \geqslant U_h(b)} p^{DM}\left(U_h\right) \geqslant \sum_{U_h \in \mathbf{\Omega}:\; U_h(b) \geqslant U_h(a)} p^{DM}(U_h),
$$
\noindent that is, if the probability mass of preference models $U_h\in \mathbf{\Omega}$ for which $U_h(a) \geqslant U_h(b)$ is not lower than the probability mass of preference models for which $U_h(b) \geqslant U_h(a)$. Of course, in this second case, to build the relation $\succsim^{DM}$, to each model $U_h\in \mathbf{\Omega}$, a mass probability $p^{DM}\left(U_h\right)$ has to be assigned. This mass probability has to be non-increasing with the distance of $U_h$ from some reference model $U_{ref}$, that is, $d^{ref}\left(U_h\right)=d\left(U_h,U_{ref}\right)$, so that, for all $U_k,U_h\in\mathbf{\Omega}$, $p^{DM}\left(U_k\right)\geqslant p^{DM}\left(U_h\right)$ iff $d\left(U_k,U_{ref}\right) \leqslant d\left(U_h,U_{ref}\right)$. Assuming the existence of a probability distribution over $\mathbf{\Omega}$ we can then consider different shapes of the DM probability distribution function $\mathbf{p}^{DM}$ shown in Table \ref{pdftab}. Note that, all probability functions we used are unimodal with a decreasing shape with respect to the distance from the reference model ${U}_{ref}$. In such distributions, the parameter $\lambda$ assumes a different meaning according to the analytic form of the distribution. For example, in the $\mathbf{p}^{Nor}$, $\lambda$ is equal to the discrete standard deviation of the distribution, while, in $\mathbf{p}^{ROC}$, it is a percentage of models, ordered with respect to the distance from the reference model, which receive a Rank Order Centroid weight (see \citealt{barron1996decision}). This means that $\mathbf{p}^{ROC}$ assigns a probability mass only to the first $\lambda$ models, $U_1,U_2,\ldots,U_{\lambda}$ in $\mathbf{\Omega}=\{U_1,U_2,\ldots,U_{|\mathbf{\Omega}|}\}$ that are renumbered in such a way that $d^{ref}\left(U_1\right)\leqslant d^{ref}\left(U_2\right)\leqslant\cdots\leqslant d^{ref}\left(U_{|\mathbf{\Omega}|}\right)$. 
\begin{table}[!h]
\caption{DM unimodal probability distribution function}\label{pdftab}
\centering\footnotesize
\resizebox{0.8\textwidth}{!}{
\begin{tabular}{ll}
\hline
\hline
Notation        &   Probability mass function  \\
\hline

${\delta}_{{U}_{ref}}$ &  ${p}_{\delta}({U}_{h},\lambda) = \left\{ 
                                                       \begin{array}{lcl} 
									                                     1,  & \mbox{if}       & {U}_{h}={U}_{ref},\\
                                                       0, & \mbox{otherwise}. & \\
                                                       \end{array}
                                                       \right.$\\

$\mathbf{p}^{Nor}$ &  ${p}_{Nor}({U}_{h},\lambda) = \frac{{e}^{-\frac{{{d}^{ref}({U}_{h})}^{2}}{2{\lambda}^{2}}}}{\displaystyle\sum_{h=1}^{|\mathbf{\Omega}|}{e}^{-\frac{{{d}^{ref}({U}_{h})}^{2}}{2{\lambda}^{2}}}}$ \\

$\mathbf{p}^{Exp}$ &  ${p}_{Exp}({U}_{h},\lambda) = \frac{{e}^{-\lambda {d}^{ref}({U}_{h})}}{\displaystyle\sum_{h=1}^{|\mathbf{\Omega}|}{e}^{-\lambda {d}^{ref}({U}_{h})}}$   \\

$\mathbf{p}^{1/D}$ &  ${p}_{1/D}({U}_{h},\lambda) = \frac{{{{d}^{ref}({U}_{h})}^{-\lambda}}}{\displaystyle\sum_{h=1}^{|\mathbf{\Omega}|}{{{d}^{ref}({U}_{h})}^{-\lambda}}}$     \\

$\mathbf{p}^{ROC}$ &  ${p}_{ROC}({U}_{h},\lambda) = \left\{ 
                                    \begin{array}{lcl} 
									                  \frac{1}{\lambda}\displaystyle\sum_{j=h}^{\lambda}\frac{1}{j}, & \mbox{for} & \;h=1,\ldots,\lambda,\\
                                    0, & \mbox{for} & \;h=\lambda+1,\ldots,|\mathbf{\Omega}|.\\
                                    \end{array}
                                    \right.$\\
\hline
\end{tabular}
}
\end{table}
\end{itemize}
 
Since the methodology we are proposing aims to ``discover" $\mathbf{p}^{DM}_{},$ it is reasonable to compare the results obtained by using the approximated probability distribution.
%in terms of goodness of fitting. In this perspective we consider the distance between the estimated probability $\mathbf{p}=\left[p(U_h)\right]_{U_h\in\mathbf{\Omega}}$ and the ``true'' DM's probability $\mathbf{p}^{DM}=\left[p^{DM}(U_h)\right]_{U_h\in\mathbf{\Omega}}$ as goodness of fitting, in the sense that the lower the distance the better the fitting. In particular, we define the distance $d(\mathbf{p},\mathbf{p}^{DM})$ between the two probability distributions as follows
%\begin{equation}\label{DistanceProbabilities}
%d(\mathbf{p},\mathbf{p}^{DM})=\sum_{U_h \in \mathbf{\Omega}} \left|p(U_h) - p^{DM}(U_h)\right|.
%\end{equation}

%%%%%%%%%%%%%%%%%%%%%%%%%%%%%%%%%%%%%%%%%%%%%%%%%%%%%%%%%%%%%%%%%%%%%%%%%%%%%%%%%%%%%%%
%%% FINE PARTE DESCRITTIVA CONSIDERANDO UNA FUNZIONE DI UTILITA' GENERICA %%%%%%%%%%%%%
%%%%%%%%%%%%%%%%%%%%%%%%%%%%%%%%%%%%%%%%%%%%%%%%%%%%%%%%%%%%%%%%%%%%%%%%%%%%%%%%%%%%%%%
%The distance between the estimated probability and the true one can be read as a measure of their similarity so that, the smaller the distance, the better the model used to estimate it. Now, considering different shapes of $\mathbf{p}^{DM}_{},$ we propose to measure the goodness of fit through a metric approach (see \citealt{rachev2011probability}) by using:
%
%\begin{itemize}
	%\item a direct comparison between the distribution functions obtained by different approaches versus the real $\mathbf{p}^{DM}_{}$ through eq. (\ref{DistanceProbabilities}). It is a discrete version of the Kantorovich metric for random variable with finite expected value (see \citealt{rachev2011probability});
We compare the $RAI$ and $PWI$ matrices obtained by using our proposals with the $RAI$ and $PWI$ matrices obtained considering $\mathbf{p}^{DM}_{}$. \\
    On the one hand, an aspect that can matter in comparing two $RAI$ matrices $RAI_1=[b^{r}_{1}(\cdot)]$ and $RAI_2=[b^r_{2}(\cdot)]$, with $b^r_{1}(\cdot)$ and $b^r_{2}(\cdot)$ representing the rank acceptability index for the $r$-th position in $RAI_1$ and $RAI_2$, respectively, is how much they differ with respect to the distributions of the first $s\leqslant m$ positions. In this sense, we consider the distance shown in eq. (\ref{avgdist})
\begin{equation}\label{avgdist}
{d}_{}(s,{RAI}_{1},{RAI}_{2}) = \frac{\displaystyle\sum^{s}_{r=1}\left(\sum_{a \in A} \left|b_1^r(a) - b_2^r(a)\right| \right)}{s}          
    \end{equation}
    \noindent where $s$ represents the number of positions taken into account. Observe that $d(s,RAI_1,RAI_2)$ can be seen as the distance induced by the ${L}_{1}$ norm \citep{guide2006infinite} in the space of the first $s$ columns of $RAI$ matrices. In the following, we shall consider $s=3.$ \\
    %An example of the computation of this measure is given in Example \ref{distraiex1}. \\% is ${\mathbf{\theta}}_{}=\int_{-\infty}^{\infty}\left(|{f}^{1}_{} - {f}^{2}_{}| \right)dx$. 
    On the other hand, given two pairwise winning index matrices $PWI_1=[p_1(a,b)]$ and $PWI_2=[p_2(a,b)]$, we can compute their distance as follows:
		\begin{equation}\label{PWIDistance}
		{d}_{}({PWI}_{1},{PWI}_{2}) = \frac{{||{PWI}_{1} - {PWI}_{2}||}_{1}}{|A|(|A|-1)}=	\frac{\displaystyle\sum_{\substack{\{a,b \}\in A\times A,\\ a\neq b}}\left|p_1(a,b)-p_2(a,b)\right|}{|A|(|A|-1)}
		\end{equation}
		\noindent being an average of the distance between the pairwise winning indices of all pairs of alternatives at hand. \\
		
		Let us conclude this section by observing that a pairwise winning indices matrix ${PWI}$ with entries being binary values represents a complete order $\succsim$ iff the following constraints are satisfied:
		$$
\begin{array}{l}
\left.
\begin{array}{l}
{p}(a,b) + {p}(b,a) = 1,\;\;\forall (a,b)\in A\times A\\[0,2cm]
{p}(a,b) + {p}(b,c) \leqslant {p}(a,c) + 1\;\;\forall a,b,c\in A.
\end{array}
\right\}
\end{array}
$$
\noindent The first constraint implies that $\succsim$ is complete and asymmetric, while, the second one implies the transitivity of $\succsim$. In this case, considering the two complete orders $\succsim_{1}$ and $\succsim_{2}$ represented by the matrices ${PWI}_{1}$ and ${PWI}_2$, respectively, we can observe that their distance computed by eq. (\ref{PWIDistance}) is equivalent to the Kendall-Tau distance between the two orders \citep{sidney1957nonparametric}. Consequently, we can say that the distance $d(PWI_1, PWI_2)$ represents an extension of the Kendall-Tau distance to probabilistic preferences.

\subsection{Computational experiments}\label{goff}%
%%%%%%%%%%%%%%%%%%%%%%%%%%%%%%%%%%%%%%%%%%%%%%%%%%%
To test our proposals we performed 1,000 runs of Algorithm \ref{SimulationAlgorithm} which steps are described in the following lines:

\begin{algorithm*}
\caption{\label{SimulationAlgorithm}}
\begin{algorithmic}
\REQUIRE Number of alternatives ($m$), number of criteria ($n$), number of pairs of alternatives to be compared by the artificial DM ($z$), artificial DM's probability distribution $\mathbf{p}^{DM}$ on $\mathbf{\Omega}$
\REPEAT
\STATE 1. Generate a performance matrix of $m$ alternatives and $n$ criteria
\STATE 2. Build the artificial DM's reference model
\STATE 3. Sample a certain number of value functions in $\mathbf{\Omega}$
\STATE 4. Compute the probability distribution $\mathbf{p}^{DM}$ on $\mathbf{\Omega}$ and the corresponding $RAI_{DM}$ and $PWI_{DM}$ matrices  
\STATE 5. Elicit the artificial DM's preferences
\STATE 6. Apply the considered methods and compute the corresponding $RAI_{Method}$ and $PWI_{Method}$ matrices
\STATE 7. Compute the distance between $RAI_{Method}$ and $RAI_{DM}$ as well as between $PWI_{Method}$ and $PWI_{DM}$
\UNTIL 1,000 runs have not been performed
\STATE 8. Compute statistics and perform the two-sample Kolmogorov-Smirnov test
\end{algorithmic}
\end{algorithm*}

\begin{description}
 \item[\textbf{1}.] build an evaluation matrix $E$ composed of $m$ rows and $n$ columns. Each entry of the matrix is sampled uniformly in the interval $\left[0,1\right]$. The $i$-th row in $E$ represents the vector of the evaluations of alternative $a_i$ on the criteria at hand. $E$ is built in such a way that all alternatives are non-dominated; 
\item[\textbf{2}.] assume that the artificial DM's value function is a weighted sum. For this reason, we randomly generate a weight vector $\mathbf{w}_{ref}=\left[w_{ref,1},\ldots,w_{ref,n}\right]$ ($\mathbf{w}_{ref}\in\left[0,1\right]^n: \displaystyle \sum_{i=1}^{n}w_{ref,i}=1$) defining the DM's reference model $U_{ref}=\mathbf{w}_{ref}$ on the basis of which the probability distribution $\mathbf{p}^{DM}$ shown in Table \ref{pdftab} will be computed;
\item[\textbf{3}.] collect a well distributed random sample $\mathbf{\Omega}$ of weight vectors $\mathbf{w}$ in the space \\ $\mathbf{W}=\left\{\mathbf{w}\in[0,1]^{n}: \displaystyle\sum_{i=1}^{n}w_i=1\right\}$; 
\item[\textbf{4}.] considering the chosen probability distribution $\mathbf{p}^{DM}$ in Table \ref{pdftab}, the rank acceptability indices $b_{DM}^{r}(a)$ for all $a\in A$ and for all $r=1,\ldots,|A|$, as well as the pairwise winning indices $p_{DM}(a,b)$ with $a,b \in A$ are computed and collected in the $RAI_{DM}$ and the $PWI_{DM}$ matrices. Analogously, the $RAI_{Uniform}$ and $PWI_{Uniform}$ matrices are computed considering a uniform distribution over $\mathbf{\Omega}$;
\item[\textbf{5}.] we extract $z$ pairs of alternatives $(a,b) \in A \times A$ on which we assume the artificial DM has to express its preferences. In particular, we select the first $z$ pairs $(a,b)\in A\times A$ that minimize $\left|p_{Uniform}(a,b)-p_{Uniform}(b,a)\right|$. In this way, we select the ``most informative" pairs since the difference $\left|p_{Uniform}(a,b)-p_{Uniform}(b,a)\right|$ is minimal when $p_{Uniform}(a,b)\cong p_{Uniform}(b,a)\cong 0.5$ and, therefore, almost half of the compatible models are in favor of the preference of $a$ over $b$, while the others are in favor of the preference of $b$ over $a$. For each selected pair $(a,b)$, the preference of the artificial DM on it is obtained by looking at its comparison in $PWI_{DM}$: if $p_{DM}(a,b)>0.5$, then, $a\succ_{P_r}b$, while $b\succ_{P_r}a$ if $p_{DM}(b,a)>0.5$. In the rare case in which $p_{DM}(a,b)=p_{DM}(b,a)=0.5$, then the preference information of the DM is $a\succsim_{P_r}b$ and $b\succsim_{P_r}a$, that is, $a$ and $b$ are indifferent;
%The following example \ref{qustpr} illustrates this procedure.
%\begin{example}\label{qustpr}
%Let us consider the ${PWI}_{Uniform}$ matrix shown in Table \ref{uniex}. The first couple we check in the pairwise comparison matrix ${PWI}_{DM}$  is $\{a,c\}$ since $\left|p_{Uniform}(a,c)-p_{Uniform}(c,a)\right|=0.1$ and, indeed, $\left|p_{Uniform}(a,c)-p_{Uniform}(c,a)\right|=\min_{a,b\in A} \left|p_{Uniform}(a,b)-p_{Uniform}(b,a)\right|$. Note that, if in the ${PWI}_{DM}$ matrix ${p}_{DM}(a,c)\geqslant 0.50$, then we assume that the DM provides the probabilistic preference $a\succsim_{Pr}c.$

%\begin{table}[!h]
%\caption{A ${PWI}_{Uniform}$ matrix}\label{uniex}
%\centering\footnotesize
%\resizebox{0.25\textwidth}{!}{
%\begin{tabular}{rcccc}
%\hline
%     & $a$     & $b$     & $c$  & $d$ \\
%\hline
%$a$  &         & $0.10$ & $0.45$ & $0.70$ \\
%$b$  & $0.90$  &        & $0.80$ & $1.00$ \\
%$c$  & $0.55$  & $0.20$ &        & $0.40$ \\
%$d$  & $0.30$  & $0.00$ & $0.60$ &        \\
%\hline
%\end{tabular}
%}
%\end{table}
%\end{example}
\item[\textbf{6}.] considering the artificial DM's preference information built in step \textbf{5.} we apply the $SSOR^{ACG}$, $SSOR^{ACG}_{pl}$ and $SSOR^{ACG}_{nl}$ methods as described in Sections \ref{mainp}-\ref{paramp}, respectively and, then, we compute the corresponding $RAI_{Method}$ and $PWI_{Method}$ matrices; in some cases, we considered also the SSOR method presented in \cite{CorrenteEtAl2016KNOSYS} and the logistic model that we will briefly recall in Section \ref{unires}. Regarding $SSOR^{ACG}_{pl}$, let us observe that the probability distribution built on $\mathbf{\Omega}$ by the method is defined by three breakpoints: 
$$d_1=0,\;\;\; d_2=\frac{\displaystyle\max_{t=1,\ldots,|\mathbf{\Omega}|}d^{ref}(U_t)}{2},\;\;\; d_3=\displaystyle\max_{t=1,\ldots,|\mathbf{\Omega}|}d^{ref}(U_t);
$$
\begin{note}
As described in detail in Section \ref{piecep} (the same holds for $SSOR^{ACG}$ described in Section \ref{mainp}) the probability distributions computed by $SSOR^{ACG}_{pl}$ is the barycenter of a sample of probability distributions compatible with constraints in $E^{ACG}_{pl}$. This implies that, solving the LP problem (\ref{lpp2}), $E^{ACG}_{pl}$ is feasible and $\varepsilon^{*}_{pl}>0$. If for the preference information provided by the artificial DM there is no compatible value function ($E^{ACG}_{pl}$ is infeasible or $\varepsilon^{*}_{pl}\leqslant 0$), then, we cannot compute the mentioned barycenter and, therefore, we consider as probability distribution the one obtained solving the LP problem (\ref{lpp2}). It does not satisfy all constraints in $E^{ACG}_{pl}$ but, however, it is the maximally discriminant one.
\end{note}
\item[\textbf{7}.] we compute the distances described in Section \ref{illex} between, on the one hand, ${RAI}_{DM}$ and $RAI_{Method}$ and, on the other hand, between ${PWI}_{DM}$ and $PWI_{Method}$;
\item[\textbf{8}.] after 1,000 runs have been performed, we show the following results: 
\begin{itemize}
\item mean and standard deviation of the distances computed in step $\mathbf{7.}$; the lower the mean and the standard deviation, the better is the considered method in approximating the artificial DM's preferences; for brevity, in the following, we shall speak about the distance between $RAI_{DM}$ and $RAI_{Method}$ and between $PWI_{DM}$ and $PWI_{Method}$ omitting the term ``matrices"; analogously, we shall say also the distance between RAI matrices and distance between PWI matrices to denote, on the one hand, the distance between ${RAI}_{DM}$ and $RAI_{Method}$ and, on the other hand, the distance between ${PWI}_{DM}$ and $PWI_{Method}$;
\item two versions of the two-sample Kolmogorov-Smirnov test \citep{Conover1999} on the distances computed in step $\mathbf{7.}$ have been performed at the 5\% significance level. In the first (``equal test"), we test the \textit{null hypothesis} of equality between two empirical distribution functions, say $F_1$ and $F_2$ versus the \textit{alternative hypothesis} that the two distributions are different. In the second (``greater test"), we test the \textit{null hypothesis} that $F_{1}$ is smaller than or equal to $F_{2}$, versus the \textit{alternative hypothesis} that $F_{1}$ is ``greater'' than $F_{2}$. Observe that in terms of cumulative distribution functions of the $1,000$ distances, $F_1$ is greater than $F_2$ if $F_1$ is more concentrated on smaller distances than $F_2$, so that, we can say that, in this case, $F_1$ fits the DM's probability distribution better than $F_2$. Consequently, in simple words, if the null hypothesis of equality between empirical distributions can be rejected in favor of the alternative hypothesis that $F_1$ is greater than $F_2$, we can conclude that statistically $F_1$ is better than $F_2$. Observe that the above-mentioned tests are performed in a sequential way. At first, the null hypothesis of equality between $F_1$ and $F_2$ is checked. If it is not rejected, then, the two distributions are equal and that's all. In the opposite case, that is, the null hypothesis is rejected (the two distributions are not equal), the alternative hypothesis is tested to check if one distribution is greater than the other. Therefore, in the tables below (and in those included in the supplementary material), the value of $h$ in correspondence of the ordered pair of distributions $\left(F_1,F_2\right)$ has to be read in the following way: 
\begin{itemize}
\item for the equal test $h=1$ means that the null hypothesis is rejected and, therefore, $F_{1}$ and $F_{2}$ are not equal, 
\item for the equal test $h=0$ means that the null hypothesis is not rejected and, therefore, $F_{1}$ and $F_{2}$ are equal,
\item $h=1$ for the greater test means that the null hypothesis is rejected and, consequently, $F_1$ is greater than $F_2$,
\item $h=0$ for the greater test means that the null hypothesis is not rejected and, consequently, $F_1$ is smaller than or equal to $F_2$.
\end{itemize}  
\item percentage of ``correct comparisons" provided by the method: to further check the capabilities of each method to discover the DM's preference model, we computed the percentage of pairs of alternatives (not included between the $z$ provided as reference examples from the DM) correctly compared by the artificial DM. Formally, considering the pair of alternatives $(a,b)\in A\times A$ on which the artificial DM did not provide any preference information (we shall briefly call a pair of this type \textit{no reference pair}), we can state that the considered method replies the comparison provided by the artificial DM on $(a,b)$ iff the probabilistic preferences of the DM ($p_{DM}(a,b)$) correspond to the probabilistic preferences of the induced probability mass comparisons ($p_{Method}(a,b)$). This means that for both probabilistic preferences $a$ is preferred to $b$ ($p_{DM}(a,b) \geqslant 0.5$ and $p_{Method}(a,b) \geqslant 0.5$) or for both probabilistic preferences $b$ is preferred to $a$ ($p_{DM}(a,b) \leqslant 0.5$ and $p_{Method}(a,b) \leqslant 0.5$). This is equivalent to the condition
\begin{equation}\label{CorrectComparison}
\left(p_{DM}(a,b)-0.5\right)\cdot\left(p_{Method}(a,b)-0.5\right)\geqslant 0. 
\end{equation}
Denoting by $z_{Method}^{Correct}$ the number of no reference pairs correctly compared by the considered method and by $\overline{z}=\binom{m}{2}-z$ the total number of no reference pairs, the percentage of correct comparisons is therefore: 
\begin{equation}\label{CorrectComparisonPercentage}
Correct^{\%}_{Method}=\frac{z_{Method}^{Correct}}{\overline{z}}.
\end{equation}  
\end{itemize}
\end{description}

%%%%%%%%%%%%%%%%%%%%%%%%%%%%%%%%%%%%%%%%%%%%%%%%%
\section{Results of the performed experiments}\label{Results}%%%
%%%%%%%%%%%%%%%%%%%%%%%%%%%%%%%%%%%%%%%%%%%%%%%%%
In this section, we shall present the results of the experiments described in detail in the previous section. In particular, in Section \ref{unires}, we shall assume that the artificial DM is evaluating alternatives by considering a single value function sampled as in step $\mathbf{2.}$ of Algorithm \ref{SimulationAlgorithm}, while, in Section \ref{mulres} we shall show the results of the experiments assuming that (i) the artificial DM is evaluating the alternatives considering not a single value function but a plurality of them and, (ii) on the set composed of these value functions, a probability distribution, to be discovered by the considered methods, is defined. As a setting for these experiments, we considered $m=8$, $n=4$, and $z=4$. In Section \ref{SensitivitySection} we present an in-depth sensitivity analysis with respect to the values of the considered parameters and the form of the artificial DM's preference model. 

%%%%%%%%%%%%%%%%%%%%%%%%%%%%%%%%%%%%%%%%%%%%%%%%%%%%%%%%%%%%%%%%%%%%%%%%
\subsection{Unique value function (${\delta}_{U_{ref}}$)}\label{unires}%
%%%%%%%%%%%%%%%%%%%%%%%%%%%%%%%%%%%%%%%%%%%%%%%%%%%%%%%%%%%%%%%%%%%%%%%%
In this section, we assume that the artificial DM evaluates the alternatives using a unique model, being the one sampled in step $\mathbf{2.}$ of Algorithm \ref{SimulationAlgorithm} and, consequently, ${\mathbf{p}}^{DM}_{}={\delta}_{{U}_{ref}}$. Here ${RAI}_{DM}$ and ${PWI}_{DM}$ are binary matrices, while, this is not the case, in general, for $RAI$ and $PWI$ matrices obtained by our three methods. Since the artificial DM evaluates using a single value function, we can look at our proposals as estimation methods via distribution fitting. For this reason, we compare the obtained results with those got by a logistic regression model (see \citealt{agresti2003categorical}) which is the standard approach in economics for such a type of problem. We shall therefore assume that, given $a,b\in A$, the probability of weak preference of $a$ over $b$ ($a\succsim b$) is the one shown in the following equation

\begin{equation}\label{logiteq}
P(a\succsim b) = \frac{1}{1+{e}^{{\beta}_{1}({g}_{1}(b)-{g}_{1}(a))+\cdots+{\beta}_{n}({g}_{n}(b)-{g}_{n}(a))}}
\end{equation}

\noindent where $\beta_{1},\ldots,\beta_{n}$ are estimated by using the preferences given by the artificial DM as described in step \textbf{5.} of Algorithm \ref{SimulationAlgorithm}. We then use the estimated model to assign the preference probability to the no reference pairs of alternatives. We refer to this model as $Logistic$.  

Looking at the results of the performed experiments in this case, we observe the following:
\begin{itemize}
\item  In Tables \ref{RAIUnique} and \ref{PWIUnique} we show the mean and standard deviation of the distances between, on the one hand, ${RAI}_{DM}$ and $RAI_{Method}$ and, on the other hand, between ${PWI}_{DM}$ and $PWI_{Method}$ (in the tables, we put in bold and italics the best (minimum) mean and standard deviation values). In both cases (here and in the following) we have that 
$$
Method\in\{Uniform, Logistic, SSOR, SSOR^{ACG}, SSOR_{pl}^{ACG}, SSOR_{nl}^{ACG}\}.
$$ 
One can see that $SSOR^{ACG}_{nl}$ has the minimal mean distance both considering the RAI and the PWI matrices. However, it presents the greatest standard deviation considering the distance between the RAI matrices, while the standard deviation of the distances between the PWI matrices is comparable with the others. \\
As to the $Logistic$ model, let us observe that, it can be used only to compute the PWI matrix since, as shown in eq. (\ref{logiteq}), it gives the probability of weak preference of an alternative over another but it cannot be used to compute the rank acceptability indices. However, let us observe that the mean distance between $PWI_{DM}$ and $PWI_{Logistic}$ is worse than the one obtained by all other methods apart from the uniform one. The model presents also a very high standard deviation compared to the other methods (almost three times worse);
\begin{table}[!h]
\begin{center}
\caption{Mean and standard deviation of the distances between $RAI_{DM}$ and $RAI_{Method}$ and between $PWI_{DM}$ and $PWI_{Method}$ for the unique model: ${\mathbf{p}}^{DM}={\delta}_{{U}_{ref}}$}\label{RAIPWIUnique}
\subtable[Distance between $RAI_{DM}$ and $RAI_{Method}$\label{RAIUnique}]{%
\resizebox{0.3\textwidth}{!}{
\begin{tabular}{lcc}
\hline
Model            & mean    & stdv    \\     
\hline
Uniform             & 1.518 &  \textbf{\textit{0.179}} \\[1mm]
${SSOR}$            & 1.415 &  0.204 \\[1mm]
${SSOR}^{ACG}$      & 1.407 &  0.223 \\[1mm]
${SSOR}_{pl}^{ACG}$ & 1.425 &  0.219 \\[1mm]
${SSOR}^{ACG}_{nl}$ & \textbf{\textit{1.026}} & 0.507 \\[1mm]
\hline
\end{tabular}
}
}
\qquad\qquad
\subtable[Distance between $PWI_{DM}$ and $PWI_{Method}$\label{PWIUnique}]{%
\resizebox{0.3\textwidth}{!}{
\begin{tabular}{lcc}
\hline
Model            & mean    & stdv    \\     
\hline
Logistic            & 0.328 & 0.206 \\[1mm]
Uniform             & 0.355 & 0.074 \\[1mm]
${SSOR}$            & 0.310 & \textbf{\textit{0.065}} \\[1mm]
${SSOR}^{ACG}$      & 0.292 & 0.071 \\[1mm]
${SSOR}_{pl}^{ACG}$ & 0.305 & 0.074 \\[1mm]
${SSOR}^{ACG}_{nl}$ & \textbf{\textit{0.147}} & 0.080 \\[1mm]
\hline
\end{tabular}
}
}
\end{center}
\end{table}

\item Table \ref{KSUnique} shows the results of the Kolmogorov-Smirnov test applied to the distances between $RAI_{DM}$ and $RAI_{Method}$ (see Tables \ref{KSRAIUniqueEqual} and \ref{KSRAIUniqueLarger}) and to the distances between $PWI_{DM}$ and $PWI_{Method}$ (see Tables \ref{KSPWIUniqueEqual} and \ref{KSPWIUniqueLarger}). 
\renewcommand\arraystretch{1.7}
 
\begin{table}[H]
\begin{center}
\caption{Kolmogorov-Smirnov test for the distances between, on the one hand, $RAI_{DM}$ and $RAI_{Method}$ matrices and, on the other hand, between $PWI_{DM}$ and $PWI_{Method}$ matrices. $\mathbf{p}^{DM}=\delta_{U_{ref}}$.\label{KSUnique}}
\subtable[Distance between RAI matrices:``Equal" Test \label{KSRAIUniqueEqual}]{%
\resizebox{0.48\textwidth}{!}{
\begin{tabular}{lcccc}
\hline
\textbf{h/p-value} & $SSOR$ & $SSOR^{ACG}$ & $SSOR^{ACG}_{pl}$ & $SSOR^{ACG}_{nl}$ \\
\hline
$Uniform$             & 1/2,96E-23 & 1/1,11E-21 & 1/1,29E-16 & 1/3,8E-138 \\[1mm]
${SSOR}$            & $\blacksquare$ & 0/0,393527 & 0/0,066631 & 1/1,6E-99 \\[1mm]
${SSOR}^{ACG}$      & $\blacksquare$ & $\blacksquare$ & 0/0,257511 & 1/1,44E-88 \\[1mm]
${SSOR}_{pl}^{ACG}$ & $\blacksquare$ & $\blacksquare$ & $\blacksquare$ & 1/5,04E-97 \\[1mm]
\hline
\end{tabular}
}
}
\;
\subtable[Distance between RAI matrices: ``Greater" test\label{KSRAIUniqueLarger}]{%
\resizebox{0.48\textwidth}{!}{
\begin{tabular}{lccccc}
\hline
\textbf{h/p-value} & $Uniform$ & $SSOR$ & $SSOR^{ACG}$ & $SSOR^{ACG}_{pl}$ & $SSOR^{ACG}_{nl}$ \\    
\hline
$Uniform$           & $\blacksquare$   & 0/1     & 0/0,998989 & 0/0,998989 & 0/0,232195 \\[1mm]
${SSOR}$            & 1/1,48E-23 & $\blacksquare$   & 0/0,531526 & 1/0,033317 & 0/0,052409 \\[1mm]
${SSOR}^{ACG}$      & 1/5,55E-22 & 0/0,198309 & $\blacksquare$   & 0/0,129033 & 0/0,052409 \\[1mm]
${SSOR}_{pl}^{ACG}$ & 1/6,47E-17 & 0/0,746591 & 0/0,995963 & $\blacksquare$   & 0/0,079818 \\[1mm]
${SSOR}_{nl}^{ACG}$ & 1/1,9E-138 & 1/8,2E-100 & 1/7,22E-89 & 1/2,52E-97 & $\blacksquare$ \\[1mm]
\hline
\end{tabular}
}
}
\subtable[Distance between PWI matrices:``Equal" Test \label{KSPWIUniqueEqual}]{%
\resizebox{0.48\textwidth}{!}{
\begin{tabular}{lccccc}
\hline
\textbf{h/p-value} & $SSOR$ & $Logistic$ & $SSOR^{ACG}$ & $SSOR^{ACG}_{pl}$ & $SSOR^{ACG}_{nl}$ \\
\hline
$Uniform$           & 1/1,12E-32 & 1/3,56E-63 & 1/2,71E-58 & 1/1,3E-37 & 1/6,9E-306 \\[1mm]
${SSOR}$            & $\blacksquare$   & 1/2,67E-50 & 1/1,98E-08 & 1/0,00498 & 1/7,1E-265 \\[1mm]
$Logistic$          & $\blacksquare$   & $\blacksquare$   & 1/1,21E-38 & 1/1,1E-39 & 1/2,9E-98 \\[1mm]
${SSOR}^{ACG}$      & $\blacksquare$   & $\blacksquare$   & $\blacksquare$   & 1/0,000664 & 1/6,1E-220 \\[1mm]
${SSOR}_{pl}^{ACG}$ & $\blacksquare$   & $\blacksquare$   & $\blacksquare$   & $\blacksquare$   & 1/1,9E-234 \\[1mm]
\hline
\end{tabular}
}
}
\;
\subtable[Distance between PWI matrices: ``Greater" test\label{KSPWIUniqueLarger}]{%
\resizebox{0.48\textwidth}{!}{
\begin{tabular}{lcccccc}
\hline
\textbf{h/p-value} & $Uniform$ & $SSOR$ & $Logistic$ & $SSOR^{ACG}$ & $SSOR^{ACG}_{pl}$ & $SSOR^{ACG}_{nl}$ \\    
\hline
$Uniform$           & $\blacksquare$   & 0/1     & 1/1,35E-21 & 0/1     & 0/1     & 0/1 \\[1mm]
${SSOR}$            & 1/5,59E-33 & $\blacksquare$   & 1/3,79E-36 & 0/0,99094 & 0/0,585712 & 0/1 \\[1mm]
$Logistic$          & 1/1,78E-63 & 1/1,33E-50 & $\blacksquare$   & 1/6,05E-39 & 1/5,48E-40 & 0/0,820209 \\[1mm]
${SSOR}^{ACG}$      & 1/1,36E-58 & 1/9,92E-09 & 1/1,17E-37 & $\blacksquare$   & 1/0,000332 & 0/1 \\[1mm]
${SSOR}_{pl}^{ACG}$ & 1/6,48E-38 & 1/0,00249 & 1/1,07E-33 & 0/0,998989 & $\blacksquare$   & 0/1 \\[1mm]
${SSOR}_{nl}^{ACG}$ & 1/3,5E-306 & 1/3,5E-265 & 1/1,45E-98 & 1/3,1E-220 & 1/9,4E-235 & $\blacksquare$ \\[1mm]
\hline
\end{tabular}
}
}
\end{center}
\end{table}

 Looking at the ``equal test" between the RAI matrices (see table \ref{KSRAIUniqueEqual}), one can observe only three pairs of methods presenting $h=0$ and for which, consequently, the corresponding distributions are equal: $\left(SSOR,SSOR^{ACG}\right)$, $\left(SSOR,SSOR^{ACG}_{pl}\right)$ and $\left(SSOR^{ACG},SSOR^{ACG}_{pl}\right)$. In all the other cases, the value $h=1$ represents a rejection of the hypothesis of equality between the considered distributions in favor of the alternative hypothesis.\\
 Considering, instead, the equal test applied to the distances between the PWI matrices (see table \ref{KSPWIUniqueEqual}), all pairs of methods are different. \\ 
 As previously described, for the cases presenting $h=1$ for the equal test, we performed the ``greater test" whose results are shown in Tables \ref{KSRAIUniqueLarger} and \ref{KSPWIUniqueLarger}. To better read the data in the tables, let us observe that the following cases can occur: 
\begin{itemize}
\item $h=1$ for the greater test applied to $(F_1,F_2)$ and $h=0$ for the greater test applied to $(F_2,F_1)$: in this case, the curve of the cumulative distribution function of $F_1$ is greater than the one corresponding to $F_2$ but it is not true the vice versa. Therefore, $F_1$ is better than $F_2$ (see Fig. \ref{h10}),
\begin{figure}[!h] 
\centering
\caption{$h=1$ for the greater test applied to $(F_1,F_2)$ and $h=0$ for the greater test applied to $(F_2,F_1)$ \label{h10}}
\includegraphics[scale=0.65]{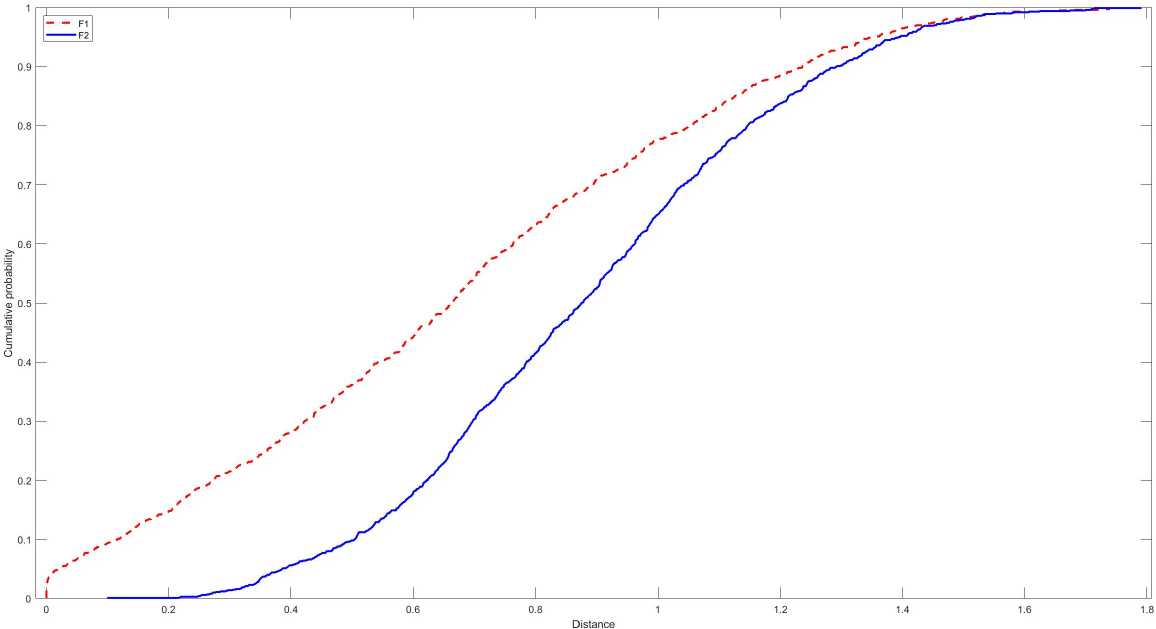}
\vspace{4mm}
\end{figure} 
\item $h=1$ for the greater test applied to $(F_1,F_2)$ as well as for the greater test applied to $(F_2,F_1)$: in this case, the curves of the cumulative distribution functions corresponding to the two methods intersect each other (see Fig. \ref{h11}) and, therefore, up to a certain distance, $F_1$ is greater than $F_2$, while, after that, $F_2$ is greater than $F_1$. In such a case, the two distributions are therefore not comparable;
%\footnote{Let us observe that the case $h=0$ for the larger test applied to the pair of methods $(F_1,F_2)$ as well as for the larger test applied to the pair of methods $(F_2,F_1)$ can never happen since this corresponds to the case in which the two distributions are equal but, in this case, we had not applied the larger test.} 
\begin{figure}[!h] 
\centering
\caption{$h=1$ both for the greater test applied to $(F_1,F_2)$ and to $(F_2,F_1)$ \label{h11}}
\includegraphics[scale=0.9]{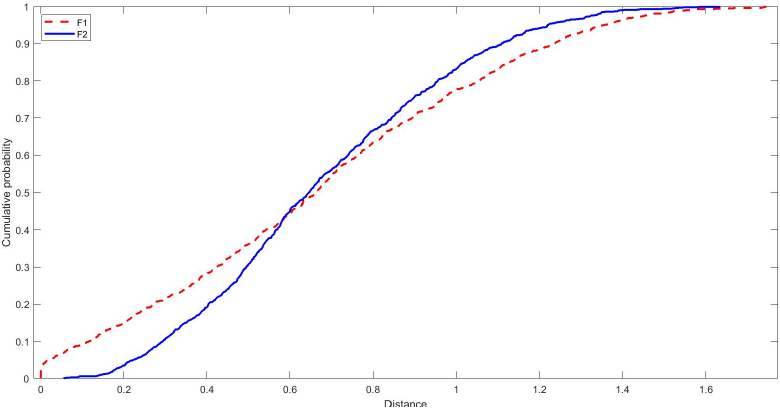}
\vspace{4mm}
\end{figure} 
\end{itemize}

Considering the greater test applied to the distances between $RAI_{DM}$ and $RAI_{Method}$ as well as between $PWI_{DM}$ and $PWI_{Method}$ one can observe that, $SSOR_{nl}^{ACG}$ is better than all other methods, while $Uniform$ is worse than all other methods. The greater test applied to the distance between the PWI matrices gives different recommendations about the three pairs of methods being equal for the equal test applied to the RAI matrices: on the one hand, $SSOR^{ACG}$ is better than $SSOR$ and $SSOR^{ACG}_{pl}$, while, on the other hand, $SSOR^{ACG}_{pl}$ is better than $SSOR$. \\
A different argumentation has to be done for the $Logistic$ model. Considering the equal and greater tests distances between $PWI_{DM}$ and $PWI_{Method}$, one can observe that its cumulative distribution function is different from all the others (see Table \ref{KSPWIUniqueEqual}). However, considering the greater test, we can say that $SSOR_{nl}^{ACG}$ is better than the $Logistic$ model but the same $Logistic$ model is incomparable with all the others.
\item To better understand the capacity of the methods to replicate the DM's preferences, in Table \ref{UniqueCorrectTable}, we report the percentage of right comparisons of no reference pairs of alternatives computed by eq. (\ref{CorrectComparisonPercentage}).
\begin{table}[!h]
\caption{Percentage of no reference pairs of alternatives correctly assigned by each method. $\mathbf{p}^{DM}=\delta_{U_{ref}}$ and $z=4$\label{UniqueCorrectTable}}
\centering\footnotesize
\resizebox{0.3\textwidth}{!}{
\begin{tabular}{lcc}
\hline
Model            & mean    & stdv       \\     
\hline
$Uniform$           & 67.98\% & 15.60\% \\
${SSOR}$            & 79.29\% & 13.28\% \\
$Logistic$          & 61.06\% & 23.61\% \\
${SSOR}^{ACG}$      & 79.65\% & 13.27\% \\
${SSOR}^{ACG}_{pl}$ & 78.48\% & 13.69\% \\
${SSOR}^{ACG}_{nl}$ & \textbf{\textit{85.23\%}} & \textbf{\textit{11.44}}\% \\
\hline
\end{tabular}
}
\end{table}
Looking at the results one can have confirmation that $SSOR_{nl}^{ACG}$ is the best among the considered methods since it is able to correctly compare $85.23\%$ of the no reference pairs of alternatives. Moreover, the method presents the lowest standard deviation meaning that the variability of the obtained results is lower than the one obtained by other methods. Looking at the worst method, instead, it seems to be the $Logistic$ since it is able to correctly compare only around 61\% of the pairs of no reference alternatives.
\end{itemize}

%%%%%%%%%%%%%%%%%%%%%%%%%%%%%%%%%%%%%%%%%%%
\subsection{Artificial DM with a probability distribution on a set of value functions}\label{mulres}%
%%%%%%%%%%%%%%%%%%%%%%%%%%%%%%%%%%%%%%%%%%%
In this section, we show the results of the experiments presented in Section \ref{illex} assuming that the artificial DM evaluates using a sample of value functions $\mathbf{\Omega}$ on which a probability distribution function $\mathbf{p}^{DM}$ among those shown in Table \ref{pdftab} is defined.

%%%%%%%%%%%%%%%%%%%%%%%%%%%%%%%%%%%%%%%%%%%%%%%%%%%%%%%%%%%%%%%%
\subsubsection{Normal distribution (${\mathbf{p}}^{Nor}$)}\label{normod}%
%%%%%%%%%%%%%%%%%%%%%%%%%%%%%%%%%%%%%%%%%%%%%%%%%%%%%%%%%%%%%%%%
In this section, we assume that the probability distribution of the artificial DM on $\mathbf{\Omega}$ is ${\mathbf{p}}^{Nor}$ (see Table \ref{pdftab}). In this case, we sample the $\lambda$ parameter in $[0.01,0.2]$ once for each of the $1,000$ runs. We obtain the following:

\begin{itemize}
\item  In Tables \ref{RAINorm} and \ref{PWINorm} we show the mean and standard deviation of the distances between, on the one hand, ${RAI}_{DM}$ and $RAI_{Method}$ and, on the other hand, between ${PWI}_{DM}$ and $PWI_{Methods}$. We can observe that $SSOR^{ACG}_{nl}$ has the lowest mean distance between both RAI and the PWI matrices. The second best method is $SSOR^{ACG}$ considering RAI and PWI mean distances, while the worst method seems to be the $Uniform$ one presenting the greatest distances. As to the standard deviation, it is very similar for all the considered methods for both distances;
\begin{table}[!h]
\begin{center}
\caption{Mean and standard deviation of the distances between $RAI_{DM}$ and $RAI_{Method}$ and between $PWI_{DM}$ and $PWI_{Method}$ for the distribution ${\mathbf{p}}^{DM}=\mathbf{p}^{Nor}$.}\label{RAIPWINorm}
\subtable[Distance between $RAI_{DM}$ and $RAI_{Method}$\label{RAINorm}]{%
\resizebox{0.3\textwidth}{!}{
\begin{tabular}{lcc}
\hline
Model            & mean    & stdv    \\     
\hline
Uniform             & 0.885 & 0.331 \\[1mm]
${SSOR}$            & 0.789 & \textbf{\textit{0.321}} \\[1mm]
${SSOR}^{ACG}$      & 0.741 & 0.342 \\[1mm]
${SSOR}_{pl}^{ACG}$ & 0.764 & 0.343 \\[1mm]
${SSOR}^{ACG}_{nl}$ & \textbf{\textit{0.684}} & 0.341 \\[1mm]
\hline
\end{tabular}
}
}
\qquad\qquad
\subtable[Distance between $PWI_{DM}$ and $PWI_{Method}$\label{PWINorm}]{%
\resizebox{0.3\textwidth}{!}{
\begin{tabular}{lcc}
\hline
Model            & mean    & stdv    \\     
\hline
Uniform             & 0.226 & 0.078 \\[1mm]
${SSOR}$            & 0.191 & 0.070 \\[1mm]
${SSOR}^{ACG}$      & 0.169 & 0.073 \\[1mm]
${SSOR}_{pl}^{ACG}$ & 0.181 & 0.073 \\[1mm]
${SSOR}^{ACG}_{nl}$ & \textbf{\textit{0.118}} & \textbf{\textit{0.056}} \\[1mm]
\hline
\end{tabular}
}
}
\end{center}
\end{table}

\item Considering the Kolmogorov-Smirnov test applied to the distances between $RAI_{DM}$ and $RAI_{Method}$ and to the distance between $PWI_{DM}$ and $PWI_{Method}$ we can observe the following\footnote{To save space we do not report here the tables containing the indicators and the $p$-values for the considered tests but they have been provided as supplementary material.}: on the one hand, the equal test applied to the distances between $RAI_{DM}$ and $RAI_{Method}$ shows that the pairs of methods $(Uniform,SSOR)$ and $(SSOR^{ACG},SSOR^{ACG}_{pl})$ are equivalent since the corresponding distributions are equal; on the other hand, considering the equal test for the distances between $PWI_{DM}$ and $PWI_{Method}$, again $(Uniform,SSOR)$ and $(SSOR^{ACG},SSOR^{ACG}_{pl})$ are equivalent. \\ 
Going to the results of the greater test applied to the distances between $RAI_{DM}$ and $RAI_{Method}$, our three proposals $SSOR^{ACG}$, $SSOR^{ACG}_{pl}$ and $SSOR^{ACG}_{nl}$ are better than $Uniform$ and $SSOR$ (being equivalent for the equal test). Considering the greater test for the distance between $PWI_{DM}$ and $PWI_{Method}$, instead, $SSOR^{ACG}_{nl}$ is better than all other methods, while $SSOR^{ACG}_{pl}$ and $SSOR^{ACG}$ are better than both $Uniform$ and $SSOR$ methods;
\item To better understand the capacity to replicate the artificial DM's preferences, in Table \ref{NormCorrectTable}, we report the percentage of right comparisons of no reference pairs of alternatives computed by eq. (\ref{CorrectComparisonPercentage}).
\renewcommand\arraystretch{1.7}
\begin{table}[!h]
\caption{Percentage of no reference pairs of alternatives correctly compared by each method considering $\mathbf{p}^{DM}=\mathbf{p}^{Nor}$\label{NormCorrectTable}}
\centering\footnotesize
\resizebox{0.3\textwidth}{!}{
\begin{tabular}{lcc}
\hline
Model            & mean    & stdv       \\     
\hline
$Uniform$           & 77.44\% & 12.92\% \\
${SSOR}$            & 86.75\% & 10.58\% \\
${SSOR}^{ACG}$      & 87.70\% & 10.27\% \\
${SSOR}^{ACG}_{pl}$ & 86.85\% & 10.77\% \\
${SSOR}^{ACG}_{nl}$ & \textbf{\textit{89.13\%}} & \textbf{\textit{9.32\%}} \\
\hline
\end{tabular}
}
\end{table}
Again we can confirm $SSOR_{nl}^{ACG}$ as the best among the considered methods having a percentage of the $89.13\%$ to correctly compare no reference pairs of alternatives. Moreover, the method presents also in this case the lowest standard deviation meaning that the variability of the obtained results is low even if similar to the one obtained by the other methods. The lowest percentage is observable for the $Uniform$ method comparing correctly $77.44\%$ of no reference pairs of alternatives, almost 10\% less than all other methods. Moreover, $Uniform$ presents the greatest standard deviation showing a great variability.
\end{itemize}

%%%%%%%%%%%%%%%%%%%%%%%%%%%%%%%%%%%%%%%%%%%%%%%%%%%%%%%%%%%%%%%%%%%
\subsubsection{Exponential distribution (${\mathbf{p}}^{Exp}$)}\label{expmod}%
%%%%%%%%%%%%%%%%%%%%%%%%%%%%%%%%%%%%%%%%%%%%%%%%%%%%%%%%%%%%%%%%%%%
If we assume that the DM evaluates the alternatives by using an exponential distribution function (${\mathbf{p}}^{DM}={\mathbf{p}}^{Exp}$ and an integer $\lambda\in[8,13]$), the results of the experiments described in Section \ref{goff} show that:
\begin{itemize}
\item  In Tables \ref{RAIExp} and \ref{PWIExp} we show the distances between $RAI_{DM}$ and $RAI_{Method}$ and between $PWI_{DM}$ and $PWI_{Method}$ obtained applying all considered methods. As one can see, $SSOR^{ACG}$ presents the lowest mean distance between the RAI matrices (see Table \ref{RAIExp}), while $SSOR^{ACG}_{nl}$ has the lowest mean distance between the PWI matrices (see Table \ref{PWIExp}). The method presenting the greatest mean distances is the $Uniform$ one. 
\begin{table}[!h]
\begin{center}
\caption{Mean and standard deviation of the distances between $RAI_{DM}$ and $RAI_{Method}$ and between $PWI_{DM}$ and $PWI_{Method}$ for the distribution ${\mathbf{p}}^{DM}=\mathbf{p}^{Exp}$.}\label{RAIPWIExp}
\subtable[Distance between $RAI_{DM}$ and $RAI_{Method}$\label{RAIExp}]{%
\resizebox{0.3\textwidth}{!}{
\begin{tabular}{lcc}
\hline
Model            & mean    & stdv    \\     
\hline
Uniform             & 0.637 & 0.207 \\[1mm]
$SSOR$              & 0.541 & 0.176 \\[1mm]
${SSOR}^{ACG}$      & \textbf{\textit{0.478}} & \textbf{\textit{0.195}} \\[1mm]
${SSOR}_{pl}^{ACG}$ & 0.484 & 0.198 \\[1mm]
${SSOR}^{ACG}_{nl}$ & 0.538 & 0.253 \\[1mm]
\hline
\end{tabular}
}
}
\qquad\qquad
\subtable[Distance between $PWI_{DM}$ and $PWI_{Method}$\label{PWIExp}]{%
\resizebox{0.3\textwidth}{!}{
\begin{tabular}{lcc}
\hline
Model            & mean    & stdv    \\     
\hline
Uniform             & 0.181 & 0.065 \\[1mm]
$SSOR$              & 0.145 & 0.052 \\[1mm]
${SSOR}^{ACG}$      & 0.124 & 0.056 \\[1mm]
${SSOR}_{pl}^{ACG}$ & 0.129 & 0.058 \\[1mm]
${SSOR}^{ACG}_{nl}$ & \textbf{\textit{0.117}} & \textbf{\textit{0.052}} \\[1mm]
\hline
\end{tabular}
}
}
\end{center}
\end{table}
\item Checking if the above comparisons are significant from the statistical point of view, the first version of the Kolmogorov-Smirnov test (equal test) applied to the distances between $RAI_{DM}$ and $RAI_{Method}$ shows that the distributions of the pairs of methods ($SSOR^{ACG}$, $SSOR^{ACG}_{pl}$), ($SSOR^{ACG}$, $SSOR^{ACG}_{nl}$) and ($SSOR^{ACG}_{pl}$, $SSOR^{ACG}_{nl}$) are equal, while, this is not the case for all other pairs of methods. In particular, the distributions corresponding to $SSOR^{ACG}$, $SSOR^{ACG}_{pl}$ and $SSOR^{ACG}_{nl}$ are greater (and, therefore, better) than the distributions corresponding to the $Uniform$ and $SSOR$ methods. \\
Considering the distributions of the distances between $PWI_{DM}$ and $PWI_{Method}$, the equal test shows that they are equal for the pair of methods ($SSOR^{ACG}$, $SSOR^{ACG}_{pl}$). Going to the greater test, on the one hand, it shows, that our proposals are better than $Uniform$ and $SSOR$ methods but, on the other hand, it proves that $SSOR^{ACG}_{nl}$ is better than both $SSOR^{ACG}$ and $SSOR^{ACG}_{pl}$; 
\item Considering the capacity of the methods to replicate the artificial DM's preferences, in Table \ref{ExpCorrectTable}, we report the percentage of right comparisons of no reference pairs of alternatives computed by eq. (\ref{CorrectComparisonPercentage}).
\renewcommand\arraystretch{1.7}
\begin{table}[!h]
\caption{Percentage of no reference pairs of alternatives correctly compared by each method considering $\mathbf{p}^{DM}=\mathbf{p}^{Exp}$\label{ExpCorrectTable}}
\centering\footnotesize
\resizebox{0.3\textwidth}{!}{
\begin{tabular}{lcc}
\hline
Model            & mean    & stdv       \\     
\hline
$Uniform$           & 77.29\% & 13.65\% \\
$SSOR$              & 87.52\% & 10.45\% \\[1mm]
${SSOR}^{ACG}$      & 88.12\% & 10.18\% \\
${SSOR}^{ACG}_{pl}$ & 87.55\% & 10.30\% \\
${SSOR}^{ACG}_{nl}$ & \textbf{\textit{89.32\%}} & \textbf{\textit{9.94\%}} \\
\hline
\end{tabular}
}
\end{table}
Looking at the results one can see that $SSOR_{nl}^{ACG}$ has the greatest percentage of correct comparisons (89.32\%). Again, the $Uniform$ is the worst among the methods at hand having the lowest percentage of right comparisons. As to the standard deviation $SSOR^{ACG}_{nl}$ has the lowest variability, while, the $Uniform$ method has the greatest one.
\end{itemize}

%%%%%%%%%%%%%%%%%%%%%%%%%%%%%%%%%%%%%%%%%%%%%%%%%%%%%%%%%%%%%%%%%%%%%%%%%
\subsubsection{1/D distribution (${\mathbf{p}}^{1/D}$)}\label{1dmod}%
%%%%%%%%%%%%%%%%%%%%%%%%%%%%%%%%%%%%%%%%%%%%%%%%%%%%%%%%%%%%%%%%%%%%%%%%%
Performing 1,000 iterations of the procedure described in Section \ref{illex} assuming that the artificial DM's probability distribution over $\mathbf{\Omega}$ is ${\mathbf{p}}^{1/D}$ (considering $\lambda$ randomly taken in $[2.5,7.5]$ for each run), we obtain the following:

\begin{itemize}
\item  In Tables \ref{RAI1D} and \ref{PWI1D} we show the mean and standard deviation of the distances between ${RAI}_{DM}$ and $RAI_{Method}$ as well as between ${PWI}_{DM}$ and $PWI_{Method}$ for each considered method. We can observe that $SSOR^{ACG}_{nl}$ has the lowest mean distances considering the RAI and the PWI matrices. However, it presents the greatest standard deviation considering the distance between the RAI matrices. The standard deviation of the distances computed between the PWI matrices is the lowest and comparable with the ones obtained by the other methods. $SSOR^{ACG}$ performs slightly better than $SSOR^{ACG}_{pl}$ for both distances even if the difference between the two methods is not very large. The $Uniform$ is the worst method among those considered also considering this probability distribution over $\mathbf{\Omega}$;
\begin{table}[!h]
\begin{center}
\caption{Mean and standard deviation of the distances between $RAI_{DM}$ and $RAI_{Method}$ and between $PWI_{DM}$ and $PWI_{Method}$ for the distribution ${\mathbf{p}}^{DM}=\mathbf{p}^{1/D}$.}\label{RAIPWI1D}
\subtable[Distance between $RAI_{DM}$ and $RAI_{Method}$\label{RAI1D}]{%
\resizebox{0.3\textwidth}{!}{
\begin{tabular}{lcc}
\hline
Model            & mean    & stdv    \\     
\hline
$Uniform$             & 0.778 & 0.244 \\[1mm]
${SSOR}$            & 0.687 & \textbf{\textit{0.234}} \\[1mm]
${SSOR}^{ACG}$      & 0.650 & 0.248 \\[1mm]
${SSOR}_{pl}^{ACG}$ & 0.655 & 0.250 \\[1mm]
${SSOR}^{ACG}_{nl}$ & \textbf{\textit{0.628}} &  0.286 \\[1mm]
\hline
\end{tabular}
}
}
\qquad\qquad
\subtable[Distance between $PWI_{DM}$ and $PWI_{Method}$\label{PWI1D}]{%
\resizebox{0.3\textwidth}{!}{
\begin{tabular}{lcc}
\hline
Model            & mean    & stdv    \\     
\hline
$Uniform$           & 0.205 & 0.068 \\[1mm]
${SSOR}$            & 0.170 & 0.058 \\[1mm]
${SSOR}^{ACG}$      & 0.149 & 0.061 \\[1mm]
${SSOR}_{pl}^{ACG}$ & 0.155 & 0.063 \\[1mm]
${SSOR}^{ACG}_{nl}$ & \textbf{\textit{0.119}} & \textbf{\textit{0.050}} \\[1mm]
\hline
\end{tabular}
}
}
\end{center}
\end{table}

\item Applying the equal test of the Kolmogorov-Smirnov to the distributions of the distances between the RAI matrices and between the PWI matrices, one gets that $SSOR^{ACG}$ and $SSOR^{ACG}_{pl}$ are equal in both cases, while, the distributions for all other pairs of methods are not equal. Therefore, we performed the greater test of the Kolmogorov-Smirnov, getting the following evidence: (i) $SSOR^{ACG}_{nl}$ is better than $Uniform,$ $SSOR^{ACG}$ and $SSOR^{ACG}_{pl}$ considering the distance between RAI matrices as well as the distance between PWI matrices; (ii) our three proposals ($SSOR^{ACG}$ and $SSOR^{ACG}_{pl}$ and $SSOR^{ACG}_{nl}$) are better than the $Uniform$ and $SSOR$ methods;
\item Looking at Table \ref{1DCorrectTable}, one has the confirmation of the goodness of the $SSOR^{ACG}_{nl}$ since it presents the greatest percentage (88.18\%) of correct comparisons of no reference pairs of alternatives computed using eq. (\ref{CorrectComparisonPercentage}). The same percentage computed for $SSOR$, $SSOR^{ACG}$ and $SSOR^{ACG}_{pl}$ is quite similar, even if the one corresponding to $SSOR^{ACG}$ is slightly greater. The $Uniform$ method confirms to be the worst among those under consideration. \\
\renewcommand\arraystretch{1.7}
\begin{table}[!h]
\caption{Percentage of no reference pairs of alternatives correctly compared by each method considering $\mathbf{p}^{DM}=\mathbf{p}^{1/D}$\label{1DCorrectTable}}
\centering\footnotesize
\resizebox{0.3\textwidth}{!}{
\begin{tabular}{lcc}
\hline
Model            & mean    & stdv       \\     
\hline
$Uniform$           & 75.98\% & 13.98\% \\
${SSOR}$            & 86.20\% & 10.91\% \\
${SSOR}^{ACG}$      & 86.48\% & 10.62\% \\
${SSOR}^{ACG}_{pl}$ & 85.37\% & 11.20\% \\
${SSOR}^{ACG}_{nl}$ & \textbf{\textit{88.18\%}} & \textbf{\textit{10.32\%}} \\
\hline
\end{tabular}
}
\end{table}
As to the variability of the results, $SSOR^{ACG}_{nl}$ has the lowest standard deviation, while the maximal is observed in correspondence of the $Uniform$ method. The other three ($SSOR$, $SSOR^{ACG}$ and $SSOR^{ACG}_{pl}$) present similar values.
\end{itemize}

%%%%%%%%%%%%%%%%%%%%%%%%%%%%%%%%%%%%%%%%%%%%%%%%%%%%%%%%%%%%%%%%%%%%%%%%%%%%%
\subsubsection{ROC distribution (${\mathbf{p}}^{ROC}$)}\label{rocmod}%
%%%%%%%%%%%%%%%%%%%%%%%%%%%%%%%%%%%%%%%%%%%%%%%%%%%%%%%%%%%%%%%%%%%%%%%%%%%%%
Considering the probability distribution ${\mathbf{p}}^{ROC}$ over $\mathbf{\Omega}$ and a random $\lambda\in\{10,11,\ldots,100\}$ for each run, we can observe the following:
\begin{itemize}
\item  In Tables \ref{RAIROC} and \ref{PWIROC} we show the distances between, on the one hand, ${RAI}_{DM}$ and $RAI_{Method}$ and, on the other hand, between ${PWI}_{DM}$ and $PWI_{Method}$. We can observe that $SSOR^{ACG}_{nl}$ has the lowest mean distance in both cases even if it presents the greatest standard deviation for the RAI distances. One aspect that is worth noting is the difference between the best and the second best method in terms of the distances considered. In particular, the difference between the mean RAI distance observed for $SSOR^{ACG}_{nl}$ and $SSOR^{ACG}$ is 0.26, while, the mean PWI distance observed for $SSOR^{ACG}_{nl}$ is almost half of the mean PWI distance obtained for $SSOR^{ACG}$. The uniform is still the worst among the compared methods in both cases;
\begin{table}[!h]
\begin{center}
\caption{Mean and standard deviation of the distances between $RAI_{DM}$ and $RAI_{Method}$ and between $PWI_{DM}$ and $PWI_{Method}$ for the distribution ${\mathbf{p}}^{DM}=\mathbf{p}^{ROC}$}\label{RAIPWIROC}
\subtable[Distance between $RAI_{DM}$ and $RAI_{Method}$\label{RAIROC}]{%
\resizebox{0.3\textwidth}{!}{
\begin{tabular}{lcc}
\hline
Model            & mean    & stdv    \\     
\hline
Uniform             & 1.070 & 0.220 \\[1mm]
${SSOR}$            & 0.974 & \textbf{\textit{0.210}} \\[1mm]
${SSOR}^{ACG}$      & 0.938 & 0.218 \\[1mm]
${SSOR}_{pl}^{ACG}$ & 0.938 & 0.218 \\[1mm]
${SSOR}^{ACG}_{nl}$ & \textbf{\textit{0.678}} & 0.307 \\[1mm]
\hline
\end{tabular}
}
}
\qquad\qquad
\subtable[Distance between $PWI_{DM}$ and $PWI_{Method}$\label{PWIROC}]{%
\resizebox{0.3\textwidth}{!}{
\begin{tabular}{lcc}
\hline
Model            & mean    & stdv    \\     
\hline
Uniform             & 0.271 & 0.062 \\[1mm]
${SSOR}$            & 0.234 & \textbf{\textit{0.053}} \\[1mm]
${SSOR}^{ACG}$      & 0.213 & 0.055 \\[1mm]
${SSOR}_{pl}^{ACG}$ & 0.218 & 0.058 \\[1mm]
${SSOR}^{ACG}_{nl}$ & \textbf{\textit{0.120}} & 0.059 \\[1mm]
\hline
\end{tabular}
}
}
\end{center}
\end{table}
\item To check for the statistical significance of the difference between the distributions of the distances obtained for each pair of methods, we performed the equal test of the Kolmogorov-Smirnov obtaining that $SSOR^{ACG}$ and $SSOR^{ACG}_{pl}$ distributions are equivalent both considering the RAI and the PWI distances. All the other pairs of methods are, instead, different. Performing, for these pairs of methods the greater test, we got recommendations similar to those shown in the $\mathbf{p}^{DM}=\mathbf{p}^{1/D}$ case: (i) the distribution corresponding to $SSOR^{ACG}_{nl}$ is greater than the one corresponding to all other methods, while, (ii) the distributions of the distances between RAI matrices and between PWI matrices observed for all methods are greater than the one obtained from the $Uniform$ method being, therefore, the worst among them;
\item To better understand the capacity of the methods to replicate the artificial DM's preferences, in Table \ref{ROCCorrectTable}, we report the percentage of right comparisons of no reference pairs of alternatives computed by eq. (\ref{CorrectComparisonPercentage}).
\renewcommand\arraystretch{1.7}
\begin{table}[!h]
\caption{Percentage of no reference pairs of alternatives correctly assigned by each method considering $\mathbf{p}^{DM}=\mathbf{p}^{ROC}$\label{ROCCorrectTable}}
\centering\footnotesize
\resizebox{0.3\textwidth}{!}{
\begin{tabular}{lcc}
\hline
Model            & mean    & stdv       \\     
\hline
$Uniform$           & 74.43\% & 13.64\% \\
${SSOR}$            & 84.65\% & 11.17\% \\
${SSOR}^{ACG}$      & 85.24\% & 10.81\% \\
${SSOR}^{ACG}_{pl}$ & 84.29\% & 10.95\% \\
${SSOR}^{ACG}_{nl}$ & \textbf{\textit{87.91\%}} & \textbf{\textit{10.68\%}} \\
\hline
\end{tabular}
}
\end{table}
Looking at the results one can observe that $SSOR_{nl}^{ACG}$ is the best among the considered methods having a percentage of right comparisons of no reference pairs of alternatives equal to $87.91\%$, being 2.67\% greater than the one observed for $SSOR^{ACG}$ (85.24\%). At the same time, the standard deviation observed for $SSOR_{nl}^{ACG}$ is the lowest (10.68\%) showing that the method presents results that are not very variable. The $Uniform$ is the worst among the considered methods also considering $\mathbf{p}^{DM}=\mathbf{p}^{ROC}$ both in terms of mean distance than in terms of the standard deviation of the same distances. 
\end{itemize}

%%%%%%%%%%%%%%%%%%%%%%%%%%%%%%%%%%%%%%%%%%%%%%
\section{Additional computational experiments}\label{SensitivitySection}%%%
%%%%%%%%%%%%%%%%%%%%%%%%%%%%%%%%%%%%%%%%%%%%%%
%%%%%%%%%%%%%%%%%%%%%%%%%%%%%%%%%%%%%%%%%%%%%%%%%%%%%%%%%%%%%%%%%%%%%%%
\subsection{Choosing a reference model differently from the barycenter}%%%
%%%%%%%%%%%%%%%%%%%%%%%%%%%%%%%%%%%%%%%%%%%%%%%%%%%%%%%%%%%%%%%%%%%%%%
In $SSOR^{ACG}$, $SSOR^{ACG}_{pl}$ and $SSOR^{ACG}_{nl}$ we aim to build a probability distribution on $\mathbf{\Omega}$ such that the probability attached to each value function in it is decreasing with the difference from a reference model $U_{ref}$. As previously described, in the three proposals, we assumed that the reference model is the $\mathbf{\Omega}$ barycenter, that is, $U_{ref}=U_{Bar}$. In this section, we shall describe a few alternative ways of selecting the reference model $U_{ref}$ showing also the results of the application of our proposals assuming the new reference model. In particular, for $SSOR^{ACG}_{pl}$ and $SSOR^{ACG}_{nl}$,  we considered the following alternative reference models $U_{ref}$: 
\begin{itemize}
\item \textit{The arithmetic mean of the simplex $\mathbf{W}$}, $U_{ArMe}$: it is a weight vector where all components have exactly the same weight, that is, $\mathbf{w}_{ArMe}=\left(w_{1}^{ArMe},\ldots,w_n^{ArMe}\right)$ such that $w_i^{ArMe}=\frac{1}{n}$, for all $i=1,\ldots,n$;
\item \textit{The most discriminant value function, $U_{MDisc}$:} it is the function obtained solving the LP problem (\ref{lpcert}) for which we get $\varepsilon^{*}.$ 
\end{itemize}
The two alternative ways of defining $U_{ref}$ presented above are not considered for $SSOR^{ACG}$ because, the equal test of the Kolmogorov-Smirnov shows that in most cases, the distributions of the distance corresponding to $SSOR^{ACG}$ and $SSOR^{ACG}_{pl}$ are equal. Since the computational effort required by $SSOR^{ACG}$ is very much larger than that one necessary for $SSOR^{ACG}_{pl}$, we applied $SSOR^{ACG}_{pl}$ only. \\
Only for the ``Parametric" method $SSOR^{ACG}_{nl}$, we have considered other two additional ways of defining the reference model: 
\begin{itemize}
\item \textit{Convex combination}, $U_{Conv}$: it is computed as a convex combination of the three reference models considered before, that is, $U_{Bar}$, $U_{ArMe}$ and $U_{MDisc}$. In particular, 
$$
U_{Conv}=\alpha\cdot U_{Bar}+\beta\cdot U_{ArMe}+\gamma\cdot U_{MDisc},
$$
where $\alpha,\beta,\gamma\geqslant 0$ and $\alpha+\beta+\gamma=1$. In this case, considering $U_{Bar}=\mathbf{w}_{Bar}=\left(w_{1}^{Bar},\ldots,w_n^{Bar}\right)$, $U_{ArMe}=\mathbf{w}_{ArMe}=\left(w_{1}^{ArMe},\ldots,w_n^{ArMe}\right)$ and $U_{MDisc}=\mathbf{w}_{Mdisc}=\left(w_{1}^{MDisc},\ldots,w_n^{MDisc}\right)$, we have that 
$U_{Conv}=\mathbf{w}_{Conv}=\left(w_{1}^{Conv},\ldots,w_n^{Conv}\right)$ such that 
$$
w_{i}^{Conv}=\alpha\cdot w_{i}^{Bar}+\beta\cdot w_{i}^{ArMe}+\gamma\cdot w_{i}^{MDisc} \;\;\mbox{for all}\;\;i=1,\ldots,n,
$$
\noindent and, consequently, for each $U_h\in\mathbf{\Omega}$, 
$$
d^{ref}(U_h)^{2}=\displaystyle\sum_{i=1}^{n}\left(w_{i}^{ref}-w_{i}^{h}\right)^2=\displaystyle\sum_{i=1}^{n}\left(w_{i}^{Conv}-w_{i}^{h}\right)^2=\displaystyle\sum_{i=1}^{n}\left(\alpha\cdot w_{i}^{Bar}+\beta\cdot w_{i}^{ArMe}+\gamma\cdot w_{i}^{MDisc}-w_{i}^{h}\right)^2
$$ 
\noindent where $\alpha$, $\beta$ and $\gamma$ are unknown. The probability distribution on $\mathbf{\Omega}$ is therefore found solving the programming problems (\ref{ppnor}) and (\ref{ppexp}) under the following constraints
$$
\left.
\begin{array}{l}
\lambda\geqslant 0, \\[0,2cm]
\alpha,\beta,\gamma\geqslant 0, \\[0,2cm]
\alpha+\beta+\gamma=1.
\end{array}
\right\}
$$
\item \textit{Unknown reference model}, $U_{Unkn}$: the unknown reference model is a weight vector $U_{Unkn}=\mathbf{w}_{Unkn}=\left(w_{1}^{Unkn},\ldots,w_{n}^{Unkn}\right)$ such that $w_i^{Unkn}\geqslant 0$, for all $i=1,\ldots,n$ and $\displaystyle\sum_{i=1}^{n}w_i^{Unkn}=1$ obtained solving the programming problems (\ref{ppnor}) and (\ref{ppexp}) under the following constraints: 
$$
\left.
\begin{array}{l}
\lambda\geqslant 0, \\[0,2cm]
w_{i}^{Unkn}\geqslant 0, \;\;\mbox{for all}\;\;i=1,\ldots,n, \\[0,2cm]
\displaystyle\sum_{i=1}^{n}w_i^{Unkn}=1.
\end{array}
\right\}
$$
In this case, the distance of each function $U_h\in\mathbf{\Omega}$, from the reference model is
$$
d^{ref}(U_h)=\sqrt{\displaystyle\sum_{i=1}^{n}\left(w_{i}^{ref}-w_{i}^{h}\right)}=\sqrt{\sum_{i=1}^{n}\left(w_{i}^{Unkn}-w_{i}^{h}\right)}.
$$ 
The probability masses $p(U_h)$ are computed by minimizing $d^{ref}(U_h)$ as shown in eqs (\ref{nlnor}) and (\ref{nlexp}) taking weights $w_i^{Unkn}, i=1,\ldots,n,$ as unknown variables. 
\end{itemize}
We considered these two reference functions only for the parametric method $SSOR^{ACG}_{nl}$ because the single components of the weight vectors $\mathbf{w}^{Conv}=\left[w^{Conv}_1, \ldots,w^{Conv}_n\right]$ and $\mathbf{w}^{Unkn}=\left[w^{Unkn}_1, \ldots,w^{Unkn}_n\right]$ are unknown variables and, therefore, also the distances $d^{ref}(U_h)$ are unknown. $SSOR^{ACG}$ and $SSOR^{ACG}_{pl}$ cannot be applied considering these two reference models since the constraints 
\begin{center}
``$p(U_k)\geqslant p({U_h})$, iff $d^{ref}(U_k)\leqslant d^{ref}(U_h)$, for all $U_h,U_k\in\mathbf{\Omega}$" 
\end{center}
present in the linear programming problems used to compute the corresponding probability distributions over $\mathbf{\Omega}$ are not linear anymore making the corresponding problems nonlinear too.   

For space reasons, we shall show the results of the experiments performed applying the considered methods but with the ``new" reference model in case $m=8$, $n=4$, $z=4$ and considering the five probability distributions shown in Table \ref{pdftab}. However, additional results obtained considering different numbers of alternatives or criteria, and different numbers of alternatives pairwise comparisons have been provided as supplementary material. Moreover, only the distance between $PWI_{DM}$ and $PWI_{Method}$ have been taken into account. \\
To underline the type of reference model used in $SSOR^{ACG}_{pl}$ or $SSOR^{ACG}_{nl}$, in the following tables in this section, we shall denote by $SSOR^{ACG}_{pl, ref}$ and $SSOR^{ACG}_{nl, ref}$ the piecewise linear method and the non-linear method described in sections \ref{piecep} and \ref{paramp}, respectively, considering as reference model $ref\in\{Bar, ArMe, MD, Conv, Unkn\}$. For example, $SSOR^{ACG}_{pl, Bar}$ will correspond to the $SSOR^{ACG}_{pl}$ method having the barycenter as a reference model (therefore, the one presented in Section \ref{piecep}), while $SSOR^{ACG}_{nl, Unkn}$ is the model presented in Section \ref{paramp} considering the unknown reference model described above. 

\renewcommand\arraystretch{1.8}
\begin{table}[!h]
\caption{Mean and standard deviation of the distances ${PWI}_{DM}$ and $PWI_{Method}$ considering the five probability distributions over $\mathbf{\Omega}$ shown in Table \ref{pdftab}\label{DistanceDifferentReferenceModels}}
\centering\footnotesize
\resizebox{1\textwidth}{!}{
\begin{tabular}{lcclcclcclcclcc}
\hline
\multicolumn{3}{c}{$\mathbf{p}^{DM}=\delta_{U_{ref}}$} & \multicolumn{3}{c}{$\mathbf{p}^{DM}=\mathbf{p}^{Nor}$} &
\multicolumn{3}{c}{$\mathbf{p}^{DM}=\mathbf{p}^{Exp}$} &
\multicolumn{3}{c}{$\mathbf{p}^{DM}=\mathbf{p}^{1/D}$} &
\multicolumn{3}{c}{$\mathbf{p}^{DM}=\mathbf{p}^{ROC}$} \\
\hline
Model            & mean    & stdv     &  Model            & mean    & stdv     & Model            & mean    & stdv     & Model            & mean    & stdv     & Model            & mean    & stdv \\ 
\hline
$Uniform$ & 0.369 & 0.070 & $Uniform$ & 0.233 & 0.079 & $Uniform$ & 0.179 & 0.063 & $Uniform$ & 0.203 & 0.065 & $Uniform$ & 0.269 & 0.062 \\
$SSOR^{ACG}_{pl,Bar}$ & 0.309 & \textbf{\textit{0.069}} & $SSOR^{ACG}_{pl,Bar}$ & 0.183 & 0.076 & $SSOR^{ACG}_{pl,Bar}$ & 0.130 & 0.057 & $SSOR^{ACG}_{pl,Bar}$ & 0.152 & 0.063 & $SSOR^{ACG}_{pl,Bar}$ & 0.219 & 0.057 \\
$SSOR^{ACG}_{pl,MD}$ & 0.294 & 0.071 & $SSOR^{ACG}_{pl,MD}$ & 0.182 & 0.071 & $SSOR^{ACG}_{pl,MD}$ & 0.133 & 0.063 & $SSOR^{ACG}_{pl,MD}$ & 0.151 & 0.056 & $SSOR^{ACG}_{pl,MD}$ & 0.216 & \textbf{\textit{0.053}}  \\
$SSOR^{ACG}_{pl,ArMe}$ & 0.386 & 0.150 & $SSOR^{ACG}_{pl,ArMe}$ & 0.249 & 0.149 & $SSOR^{ACG}_{pl,ArMe}$ & 0.212 & 0.278 & $SSOR^{ACG}_{pl,ArMe}$ & 0.229 & 0.274 & $SSOR^{ACG}_{pl,ArMe}$ & 0.295 & 0.205 \\
$SSOR^{ACG}_{nl,Bar}$  & 0.167 & 0.079 & $SSOR^{ACG}_{nl,Bar}$  & \textbf{\textit{0.122}} & \textbf{\textit{0.055}} & $SSOR^{ACG}_{nl,Bar}$  & \textbf{\textit{0.121}} & \textbf{\textit{0.050}} & $SSOR^{ACG}_{nl,Bar}$  & \textbf{\textit{0.120}} & \textbf{\textit{0.052}} & $SSOR^{ACG}_{nl,Bar}$  & \textbf{\textit{0.116}} & 0.058 \\
$SSOR^{ACG}_{nl,MD}$   & \textbf{\textit{0.150}} & 0.105 & $SSOR^{ACG}_{nl,MD}$ & 0.180 & 0.092 & $SSOR^{ACG}_{nl,MD}$ & 0.203 & 0.094 & $SSOR^{ACG}_{nl,MD}$ & 0.192 & 0.091 & $SSOR^{ACG}_{nl,MD}$ & 0.167 & 0.094   \\
$SSOR^{ACG}_{nl,ArMe}$ & 0.346 & 0.084 & $SSOR^{ACG}_{nl,ArMe}$ & 0.208 & 0.084 & $SSOR^{ACG}_{nl,ArMe}$ & 0.164 & 0.070 & $SSOR^{ACG}_{nl,ArMe}$ & 0.180 & 0.073 & $SSOR^{ACG}_{nl,ArMe}$ & 0.239 & 0.076  \\
$SSOR^{ACG}_{nl,conv}$ & 0.366 & 0.069 & $SSOR^{ACG}_{nl,conv}$ & 0.231 & 0.079 & $SSOR^{ACG}_{nl,conv}$ & 0.178 & 0.062 & $SSOR^{ACG}_{nl,conv}$ & 0.201 & 0.065 & $SSOR^{ACG}_{nl}$ & 0.267 & 0.062  \\
$SSOR^{ACG}_{nl,Unkn}$ & 0.368 & 0.070 & $SSOR^{ACG}_{nl,Unkn}$ & 0.232 & 0.079 & $SSOR^{ACG}_{nl,Unkn}$ & 0.179 & 0.063 & $SSOR^{ACG}_{nl,Unkn}$ & 0.202 & 0.065 & $SSOR^{ACG}_{nl,Unkn}$ & 0.268 & 0.062  \\
\hline
\end{tabular}
}
\end{table}

In Table \ref{DistanceDifferentReferenceModels} we show the mean and the standard deviation of the distance between the PWI matrices considering the ten methods obtained by changing the reference model and the five different probability distributions on the space of the sampled value functions $\mathbf{\Omega}$ shown in Table \ref{pdftab}. As one can see, apart from the unique distribution $\mathbf{p}^{DM}=\delta_{U_{ref}}$ the nonlinear method considering the barycenter as reference model presents the lowest mean distance between the PWI matrices and, therefore, is the best among the considered methods. When $\mathbf{p}^{DM}=\delta_{U_{ref}}$, the best method is $SSOR^{ACG}_{nl, MD}$, that is, the nonlinear model considering the most discriminant function as reference model even if the difference between the two mean values is not very great (0.017). Nevertheless, performing the greater test of the Kolmogorov-Smirnov over the distributions corresponding to these two methods, we get that the distribution corresponding to $SSOR^{ACG}_{nl, MD}$ is greater than the one corresponding to $SSOR^{ACG}_{nl, Bar}$ and, therefore, it is better. Going more in-depth on the results obtained by the $SSOR^{ACG}_{pl}$ and $SSOR^{ACG}_{nl}$ choosing differently $U_{ref}$, we can observe the following: 
\begin{itemize}
\item considering the $SSOR^{ACG}_{pl}$, the best choice for the reference model is the barycenter or the most discriminant function for all considered distributions. In particular, $SSOR^{ACG}_{pl,MD}$ is better than $SSOR^{ACG}_{pl,Bar}$ for all distributions apart from $\mathbf{p}^{Exp}$. However, the difference between the mean PWI distances is not great in all considered cases (the maximum mean difference is 0.0151 and it is obtained when $\mathbf{p}^{DM}=\delta_{U_{ref}}$). Choosing the arithmetic mean of the sampled value functions as a reference model is the worst among the three considered options for all probability distributions over $\mathbf{\Omega}$;
\item with respect to $SSOR^{ACG}_{nl}$, as observed above, $SSOR^{ACG}_{nl, Bar}$ is the best method for four out of the five distributions and, therefore, it is also the best among the variants considered. However, while $SSOR^{ACG}_{nl, MD}$ is slightly better than $SSOR^{ACG}_{nl, Bar}$ when $\mathbf{DM}=\delta_{U_{ref}}$, the difference between the mean distance obtained for $SSOR^{ACG}_{nl,MD}$ and for $SSOR^{ACG}_{nl,Bar}$ is quite great when $SSOR^{ACG}_{nl,Bar}$ is the best among the two. For this particular model, the two versions considering the reference model as unknown are always giving the worst results in terms of the mean distance between the PWI matrices. As to the version of the method choosing the arithmetic mean as the reference model, it is always worse than $SSOR^{ACG}_{nl, Bar}$ but, in two cases ($\mathbf{p}^{Exp}$ and $\mathbf{p}^{1/D}$), it is better than $SSOR^{ACG}_{nl, MD}$ and, therefore, it is the second best option between the $SSOR^{ACG}_{nl}$ variants for these two cases;
\end{itemize}

\renewcommand\arraystretch{1.7}
\begin{table}[!h]
\caption{Percentage of correct comparisons of no reference pairs of alternatives considering the five probability distributions over $\mathbf{\Omega}$ shown in Table \ref{pdftab}\label{PercentageDifferentReferenceModels}}
\centering\footnotesize
\resizebox{1\textwidth}{!}{
\begin{tabular}{lcclcclcclcclcc}
\hline
\multicolumn{3}{c}{$\mathbf{p}^{DM}=\delta_{U_{ref}}$} & \multicolumn{3}{c}{$\mathbf{p}^{DM}=\mathbf{p}^{Nor}$} &
\multicolumn{3}{c}{$\mathbf{p}^{DM}=\mathbf{p}^{Exp}$} &
\multicolumn{3}{c}{$\mathbf{p}^{DM}=\mathbf{p}^{1/D}$} &
\multicolumn{3}{c}{$\mathbf{p}^{DM}=\mathbf{p}^{ROC}$} \\
\hline
Model            & mean    & stdv     &  Model            & mean    & stdv     & Model            & mean    & stdv     & Model            & mean    & stdv     & Model            & mean    & stdv \\ 
\hline
$Uniform$ & 67.05\% & 15.28\% & $Uniform$ & 76.30\% & 12.99\% & $Uniform$ & 77.59\% & 12.98\% & $Uniform$ & 75.70\% & 13.38\% & $Uniform$ & 74.71\% & 13.58\% \\
$SSOR^{ACG}_{pl,Bar}$ & 77.72\% & 13.35\% & $SSOR^{ACG}_{pl,Bar}$ & 85.66\% & 11.66\% & $SSOR^{ACG}_{pl,Bar}$ & 86.39\% & 11.75\% & $SSOR^{ACG}_{pl,Bar}$ & 84.73\% & 11.26\% & $SSOR^{ACG}_{pl,Bar}$ & 83.87\% & 11.70\% \\
$SSOR^{ACG}_{pl,MD}$ & 80.37\% & 14.20\% & $SSOR^{ACG}_{pl,MD}$ & 85.44\% & 12.17\% & $SSOR^{ACG}_{pl,MD}$ & 86.12\% & 12.03\% & $SSOR^{ACG}_{pl,MD}$ & 84.47\% & 11.78\% & $SSOR^{ACG}_{pl,MD}$ & 84.14\% & 11.66\%  \\
$SSOR^{ACG}_{pl,ArMe}$ & 59.52\% & 21.35\% & $SSOR^{ACG}_{pl,ArMe}$ & 68.23\% & 21.57\% & $SSOR^{ACG}_{pl,ArMe}$ & 70.74\% & 21.26\% & $SSOR^{ACG}_{pl,ArMe}$ & 68.60\% & 20.75\% & $SSOR^{ACG}_{pl,ArMe}$ & 65.46\% & 22.20\% \\
$SSOR^{ACG}_{nl,Bar}$  & \textbf{\textit{84.72\%}} & \textbf{\textit{11.63\%}} & $SSOR^{ACG}_{nl,Bar}$  & \textbf{\textit{88.37\%}} & \textbf{\textit{10.82\%}} & $SSOR^{ACG}_{nl,Bar}$  & \textbf{\textit{88.98\%}} & \textbf{\textit{10.81\%}} & $SSOR^{ACG}_{nl,Bar}$  & \textbf{\textit{88.58\%}} & \textbf{\textit{9.84\%}}  & $SSOR^{ACG}_{nl,Bar}$  & \textbf{\textit{88.39\%}} & \textbf{\textit{10.14\%}} \\
$SSOR^{ACG}_{nl,MD}$   & 83.40\% & 16.81\% & $SSOR^{ACG}_{nl,MD}$ & 80.62\% & 15.45\% & $SSOR^{ACG}_{nl,MD}$ & 79.73\% & 17.49\% & $SSOR^{ACG}_{nl,MD}$ & 79.03\% & 18.50\% & $SSOR^{ACG}_{nl,MD}$ & 80.65\% & 15.44\%   \\
$SSOR^{ACG}_{nl,ArMe}$ & 66.23\% & 15.29\% & $SSOR^{ACG}_{nl,ArMe}$ & 76.52\% & 12.93\% & $SSOR^{ACG}_{nl,ArMe}$ & 77.63\% & 13.26\% & $SSOR^{ACG}_{nl,ArMe}$ & 75.59\% & 13.50\% & $SSOR^{ACG}_{nl,ArMe}$ & 74.57\% & 13.47\%   \\
$SSOR^{ACG}_{nl,conv}$ & 67.05\% & 15.11\% & $SSOR^{ACG}_{nl,conv}$ & 76.40\% & 12.94\% & $SSOR^{ACG}_{nl,conv}$ & 77.57\% & 12.86\% & $SSOR^{ACG}_{nl,conv}$ & 75.52\% & 13.30\% & $SSOR^{ACG}_{nl}$ & 76.64\% & 13.55\%  \\
$SSOR^{ACG}_{nl,Unkn}$ & 66.67\% & 15.16\% & $SSOR^{ACG}_{nl,Unkn}$ & 76.12\% & 13.01\% & $SSOR^{ACG}_{nl,Unkn}$ & 77.22\% & 12.97\% & $SSOR^{ACG}_{nl,Unkn}$ & 74.96\% & 13.43\% & $SSOR^{ACG}_{nl,Unkn}$ & 74.28\% & 13.53\%  \\
\hline
\end{tabular}
}
\end{table}

To check for the capability of the methods to replicate the artificial DM's preferences, we computed the percentage of no reference pairs correctly compared by each method. Looking at Table \ref{PercentageDifferentReferenceModels} one can observe that $SSOR^{ACG}_{nl, Bar}$ is undoubtedly the best among the methods at hand since it presents the greatest percentage of correct comparisons in all of the three considered probability distributions. As to the second best, the situation is slightly different since there is not a single method being the second best for all cases. In particular, $SSOR^{ACG}_{nl,MD}$ is the second best considering $\delta_{U_{ref}}$, $SSOR^{ACG}_{pl,Bar}$ is the second best considering $\mathbf{p}^{Nor}$, $\mathbf{p}^{Exp}$ and $\mathbf{p}^{1/D}$, while $SSOR^{ACG}_{pl,MD}$ is the second best considering $\mathbf{p}^{ROC}$. The results show therefore, that choosing the barycenter as reference model is the best option considering all distributions for the $SSOR^{ACG}_{nl}$ method, while it is the best option for three out of the five distributions for the $SSOR^{ACG}_{pl}$ one. As to the worst between the methods at hand, there is no doubt that it is $SSOR^{ACG}_{pl, ArMe}$ since it presents the lowest percentage of right comparisons of no reference pairs of alternatives for all considered cases. 

%%%%%%%%%%%%%%%%%%%%%%%%%%%%%%%%%%%%%%%%%%%%%%%%%%%%%%%%%%%%%%%%%%%
\subsection{Sensitivity analysis with respect to the barycenter}%%%
%%%%%%%%%%%%%%%%%%%%%%%%%%%%%%%%%%%%%%%%%%%%%%%%%%%%%%%%%%%%%%%%%%%
As observed in the previous section, $SSOR^{ACG}_{pl,Bar}$ and $SSOR^{ACG}_{nl,Bar}$, that is the versions of the $SSOR^{ACG}_{pl}$ and $SSOR^{ACG}_{nl}$ methods considering the barycenter as the reference model, perform very well with respect to the other methods. In this section, we shall perform a sensitivity analysis on the barycenter itself, that is, we shall check how the results will change if the reference model is a ``neighbor" of the barycenter. We shall then apply $SSOR^{ACG}_{pl,Bar}$ and $SSOR^{ACG}_{nl,Bar}$ considering as reference model a value function close to the barycenter and, then, we shall see how the distance between the PWI matrices, as well as the percentage of the right comparisons of no reference pairs of alternatives, will change. From a technical point of view, we implemented the procedure articulated in the following steps: 
\begin{description}
\item[Step 1)] Find the barycenter $\mathbf{w}_{Bar}=\left(w_{1}^{Bar},\ldots,w_n^{Bar}\right)$ of the set $\mathbf{\Omega}$ as described in Section \ref{mainp};
\item[Step 2)] Define the set $\mathbf{W}_{NeighBar}$ composed of the weight vectors in a ``neighbor" of $\mathbf{w}_{Bar}$, that is, the set composed of the weight vectors $\left(\overline{w}_1,\ldots,\overline{w}_n\right)$ satisfying the following constraints
$$
\left.
\begin{array}{l}
w_i^{Bar}-0.05\leqslant\overline{w}_i\leqslant{w}_{i}^{Bar}+0.05,\;\;\mbox{for all}\;\;i=1,\ldots,n,\\[0,2cm]
\displaystyle\sum_{i=1}^{n}\overline{w}_i=1,\\[0,2cm]
\overline{w}_i\geqslant 0,\;\;\mbox{for all}\;\;i=1,\ldots,n;\\[0,2cm]
\end{array}
\right\}
$$
\item[Step 3)] Sample 100 weight vectors from $\mathbf{W}_{NeighBar}$ and select the nine closest to $\mathbf{w}_{Bar}$, that is, the nine weight vectors in $\mathbf{W}_{NeighBar}$ presenting the lowest Euclidean distance from $\mathbf{w}_{Bar}$.
\item[Step 4)] Apply the procedure described in Section \ref{goff} considering $m=8$, $n=4$, $z=4$ and assuming that the artificial DM evaluates the considered alternatives using a plurality of value functions on which the five probability distributions in Table \ref{pdftab} are defined.
\end{description}
In the following, we shall present the results obtained considering $\mathbf{p}^{Exp}$ as the probability distribution of the artificial DM. The data for the other four probability distributions are provided in the supplementary material. Also in this case, we compare the results obtained by $SSOR^{ACG}_{pl}$ applying $SSOR^{ACG}_{nl}$ only. \\
In Table \ref{SensitivityBarycenter} methods $SSOR^{ACG}_{pl,Bar_k}$ and $SSOR^{ACG}_{nl,Bar_k}$, with $k\in\{1,\ldots,9\}$ represents the methods $SSOR^{ACG}_{pl}$ and $SSOR^{ACG}_{nl}$ obtained considering as reference model the $k$-th value function closest to the barycenter in $\mathbf{W}_{NeighBar}$.

\begin{table}[!h]
\begin{center}
\caption{Mean distance between the PWI matrices ($\overline{d}(PWI_{DM},PWI_{Method})$ and mean percentage of correct comparisons of no reference pairs of alternatives ($\overline{Correct}^{\%}_{Method}$) considering the distribution ${\mathbf{p}}^{DM}=\mathbf{p}^{Exp}$}\label{SensitivityBarycenter}
\subtable[$Method=SSOR^{ACG}_{pl,ref}$\label{SSORplSensitivity}]{%
\resizebox{0.4\textwidth}{!}{
\begin{tabular}{lcc}
\hline
Model            & $\overline{d}(PWI_{DM},PWI_{Method})$    & $\overline{Correct}^{\%}_{Method}$    \\     
\hline
${SSOR}^{ACG}_{pl,Bar}$          & 0.1313 & 86.15\% \\[1mm]
${SSOR}^{ACG}_{pl,Bar_1}$        & 0.1310 & 86.15\% \\[1mm]
${SSOR}^{ACG}_{pl,Bar_2}$        & 0.1317 & 86.08\% \\[1mm]
${SSOR}^{ACG}_{pl,Bar_3}$        & 0.1314 & 86.10\% \\[1mm]
${SSOR}^{ACG}_{pl,Bar_4}$        & 0.1318 & 86.03\% \\[1mm]
${SSOR}^{ACG}_{pl,Bar_5}$        & 0.1317 & 86.01\% \\[1mm]
${SSOR}^{ACG}_{pl,Bar_6}$        & 0.1321 & 85.91\% \\[1mm]
${SSOR}^{ACG}_{pl,Bar_7}$        & 0.1327 & 85.90\% \\[1mm]
${SSOR}^{ACG}_{pl,Bar_8}$        & 0.1322 & 85.91\% \\[1mm]
${SSOR}^{ACG}_{pl,Bar_9}$        & 0.1322 & 85.97\% \\[1mm]
\hline
\end{tabular}
}
}
\qquad\qquad
\subtable[$Method=SSOR^{ACG}_{nl,ref}$\label{SSORnlSensitivity}]{%
\resizebox{0.4\textwidth}{!}{
\begin{tabular}{lcc}
\hline
Model            & $\overline{d}(PWI_{DM},PWI_{Method})$    & $\overline{Correct}^{\%}_{Method}$    \\     
\hline
${SSOR}^{ACG}_{nl,Bar}$          & 0.1263 & 88.70\% \\[1mm]
${SSOR}^{ACG}_{nl,Bar_1}$        & 0.1268 & 88.17\% \\[1mm]
${SSOR}^{ACG}_{nl,Bar_2}$        & 0.1275 & 88.48\% \\[1mm]
${SSOR}^{ACG}_{nl,Bar_3}$        & 0.1280 & 87.50\% \\[1mm]
${SSOR}^{ACG}_{nl,Bar_4}$        & 0.1276 & 88.38\% \\[1mm]
${SSOR}^{ACG}_{nl,Bar_5}$        & 0.1284 & 88.02\% \\[1mm]
${SSOR}^{ACG}_{nl,Bar_6}$        & 0.1277 & 88.22\% \\[1mm]
${SSOR}^{ACG}_{nl,Bar_7}$        & 0.1271 & 88.20\% \\[1mm]
${SSOR}^{ACG}_{nl,Bar_8}$        & 0.1278 & 88.17\% \\[1mm]
${SSOR}^{ACG}_{nl,Bar_9}$        & 0.1273 & 88.03\% \\[1mm]
\hline
\end{tabular}
}
}
\end{center}
\end{table}

Looking at the data presented in Table \ref{SensitivityBarycenter} one can see that choosing as reference model a value function ``close" to the barycenter, does not affect the obtained results. Indeed, as can be seen in Table \ref{SSORplSensitivity}, for the $SSOR^{ACG}_{pl}$ model, the distance between the PWI matrices considering the barycenter as reference model is 0.1313, while, the same distance considering the nine sampled weight vectors closest to the barycenter varies between 0.1310 and 0.1327. The same can be observed for the mean percentage of correct comparisons of no reference pairs of alternatives being equal to 86.15\% for the $SSOR^{ACG}_{pl,Bar}$, while, it varies in the interval $\left[85.90\%,86.15\%\right]$ for the other nine models using the closest weight vectors to the barycenter as reference models. 

Analogous considerations can be done for the $SSOR^{ACG}_{nl}$ variants obtained by choosing as reference models the nine sampled weight vectors closest to the barycenter. On the one hand, for $SSOR^{ACG}_{nl, Bar}$, the mean distance between the PWI matrices is 0.1263, while the same distance obtained considering as reference model the nine sampled weight vectors closest to the barycenter varies in $\left[0.1268,0.1284\right]$; (ii) the mean percentage of correct comparisons of no reference pairs of alternatives for $SSOR^{ACG}_{nl,Bar}$ is 88.70\%, while the same percentage varies between 87.50\% and 88.70\% considering the other nine weight vectors. 

Similar considerations can be made taking into account the other four probability distributions $\mathbf{p}^{DM}$ shown in Table \ref{pdftab} whose results are provided as supplementary material. Both the mean distance between the PWI matrices and the mean percentage of correct comparisons of no reference pairs of alternatives are not very sensitive to the choice of a weight vector close to the barycenter as reference model.

%%%%%%%%%%%%%%%%%%%%%%%%%%%%%%%%%%%%%%%%%%%%%%%%%%%%%%%%%%%%%%%%%%%%%
\subsection{Sensitivity analysis on the number of alternatives and criteria}
%%%%%%%%%%%%%%%%%%%%%%%%%%%%%%%%%%%%%%%%%%%%%%%%%%%%%%%%%%%%%%%%%%%%%
In this section, we shall study how the results obtained by the performed experiments are dependent on the number of alternatives and criteria taken into account. In particular, we shall analyze three different cases: (5,3), (8,4), and (12,6) where the first element in the pair is the number of alternatives, while, the second element is the number of criteria. To this aim, we shall assume that the artificial DM provides four pairwise comparisons of reference alternatives ($z=4$) and that it evaluates the alternatives using the value functions in $\mathbf{\Omega}$ on which the five probability distributions shown in Table \ref{pdftab} are defined. Since we have already shown that choosing the barycenter as the reference model is the best option, we shall consider only the methods in which the reference model is, indeed, the barycenter. Moreover, for this set of experiments, shown as described in Section \ref{goff}, we considered also $Uniform$ and $SSOR$. The $Logistic$ model has been considered (as before) only in case $\mathbf{p}^{DM}=\delta_{U_{ref}}$.

\renewcommand\arraystretch{1.7}
\begin{table}[!h]
\caption{Mean distance between the PWI matrices ($\overline{d}(PWI_{DM},PWI_{Method})$ and mean percentage of correct comparisons of no reference pairs of alternatives ($\overline{Correct}^{\%}_{Method}$) considering the distribution ${\mathbf{p}}^{DM}=\delta_{U_{ref}}$\label{SensitivityConfigurationsUnique}}
\centering\footnotesize
\resizebox{1\textwidth}{!}{
\begin{tabular}{lcclcclcc}
\hline
\multicolumn{3}{c}{$(5,3)$} & \multicolumn{3}{c}{$(8,4)$} &
\multicolumn{3}{c}{$(12,6)$} \\
\hline
Model            & $\overline{d}(PWI_{DM},PWI_{Method})$    & $\overline{Correct}^{\%}_{Method}$ & Model            & $\overline{d}(PWI_{DM},PWI_{Method})$    & $\overline{Correct}^{\%}_{Method}$ & Model            & $\overline{d}(PWI_{DM},PWI_{Method})$    & $\overline{Correct}^{\%}_{Method}$ \\ 
\hline
$Uniform$ & 0.379 & 50.57\% & $Uniform$ & 0.355 & 67.98\% & $Uniform$ & 0.356 & 70.92\% \\ 
$SSOR$ & 0.299 & 80.55\% & $SSOR$ & 0.310 & 79.29\% & $SSOR$ & 0.331 & 76.63\%  \\ 
$Logistic$ & 0.191 & 68.12\% & $Logistic$ & 0.328 & 61.06\% & $Logistic$ & 0.353 & 62.39\%  \\ 
$SSOR^{ACG}$ & 0.275 & 81.93\% & $SSOR^{ACG}$ & 0.292 & 79.65\% & $SSOR^{ACG}$ & 0.320 & 76.20\%  \\ 
$SSOR^{ACG}_{pl}$ & 0.282 & 80.30\% & $SSOR^{ACG}_{pl}$ & 0.305 & 78.48\% & $SSOR^{ACG}_{pl}$ & 0.331 & 75.30\%  \\ 
$SSOR^{ACG}_{nl}$ & \textbf{\textit{0.069}} & \textbf{\textit{83.00\%}} & $SSOR^{ACG}_{nl}$ & \textbf{\textit{0.147}} & \textbf{\textit{85.23\%}} & $SSOR^{ACG}_{nl}$ & \textbf{\textit{0.276}} & \textbf{\textit{79.23\%}} \\
\hline
\end{tabular}
}
\end{table}

Looking at Table \ref{SensitivityConfigurationsUnique} one can observe that, while the mean distance between the PWI matrices increases with the ``dimension" of the configuration for all the considered methods, the same can not be said for the mean percentage of correct comparisons of no reference pairs of alternatives. For example, $Uniform$ presents a mean percentage of correct pairwise comparisons increasing with the dimension of the problem since it passes from 50.57\% of the (5,3) configuration to 67.98\% of the (8,4) configuration to the 70.92\% of the (12,6) configuration, while, a completely opposite trend can be observed for $SSOR$ presenting the greatest mean percentage for the (5,3) configuration (80.55\%) and the lowest for the (12,6) configuration (76.63\%). In other cases, instead, the mean percentage is not monotonic with respect to the dimension of the configuration. For example, for $SSOR^{ACG}_{nl}$, this percentage increases passing from $83.00\%$ of the (5,3) configuration to 85.23\% of the (8,4) configuration but, then, it decreases to 79.23\% passing to the (12,6) configuration. \\
Independently of these trends, what is more important to underline is that $SSOR^{ACG}_{nl}$ presents the lowest mean distance between the PWI matrices and the greatest mean percentage of correct comparisons of no reference pairs of alternatives in the three considered configurations. At the same time, the $Uniform$ is the method presenting the worst values both in terms of the mean distance between the PWI matrices and the mean percentage in the three analyzed configurations. \\
Similar considerations can be made taking into account the other four probability distributions (see the corresponding tables provided in the supplementary material). In particular: (i) $SSOR^{ACG}_{nl}$ has the lowest mean distance between the PWI matrices in all considered cases (so all probability distributions and all configurations) apart from the configuration (12,6) with probability distribution $\mathbf{p}^{ROC}$ for which the lowest mean distance is obtained by $SSOR$; (ii) as to the mean percentage of correct comparisons provided by the DM, we observe that the best methods are $SSOR^{ACG}_{nl}$ and $SSOR^{ACG}$. In particular: 
\begin{itemize}
\item considering $\mathbf{p}^{Nor}$, $\mathbf{p}^{1/D}$ and $\mathbf{p}^{ROC},$ $SSOR^{ACG}_{nl}$ presents the greatest mean percentage of correct comparisons of no reference pairs of alternatives for configurations (8,4) and (12,6), while, $SSOR^{ACG}$ has the greatest percentage value for configuration (5,3);
\item considering $\mathbf{p}^{Exp},$ $SSOR^{ACG}_{nl}$ presents the greatest mean percentage of correct comparisons of no reference pairs of alternatives for configurations (5,3) and (8,4), while, $SSOR$ has the greatest percentage value for configuration (12,6);
\item in cases for which $SSOR^{ACG}_{nl}$ is not the best method, the difference between the mean percentage of correct comparisons it provides and the one provided by the best method ($SSOR^{ACG}$), is always very small (lower than 3\%) meaning that the results obtained by $SSOR^{ACG}_{nl}$ in terms of the capability of discovering artificial DM's preferences are anyway quite good. 
\end{itemize}

%%%%%%%%%%%%%%%%%%%%%%%%%%%%%%%%%%%%%%%%%%%%%%%%%%%%%%%%%%%%%%%%%%%%%%%%%%
\subsection{Sensitivity analysis on the number of pairwise comparisons of reference alternatives}\label{Sensitivityz}
%%%%%%%%%%%%%%%%%%%%%%%%%%%%%%%%%%%%%%%%%%%%%%%%%%%%%%%%%%%%%%%%%%%%%%%%%%
In this section, we shall present the results obtained from a sensitivity analysis on the number of pairwise comparisons provided by the DM. The analysis is conducted considering the (8,4) configuration and the probability distributions $\mathbf{p}^{DM}=\delta_{U_{ref}}$ and $\mathbf{p}^{DM}=\mathbf{p}^{Nor}$. 

\begin{table}[!h]
\begin{center}
\caption{Mean distance between the PWI matrices ($\overline{d}(PWI_{DM},PWI_{Method})$ and mean percentage of correct comparisons of no reference pairs of alternatives ($\overline{Correct}^{\%}_{Method}$) considering the distribution ${\mathbf{p}}^{DM}=\delta_{U_{ref}}$}\label{SensitivityComparisonsUnique}
\subtable[$\overline{d}(PWI_{DM},PWI_{Method})$\label{DistanceSensitivityComparisonsUnique}]{%
\resizebox{0.4\textwidth}{!}{
\begin{tabular}{lccccc}
\hline
Model            & $z=4$ & $z=9$ & $z=14$ & $z=19$ & $z=25$    \\     
\hline
$Uniform$  & 0.355 & 0.365 & 0.374 & 0.367 & 0.372 \\
$Logistic$ & 0.328 & 0.158 & 0.075 & \textbf{\textit{0.030}} & \textbf{\textit{0.006}} \\ 
$SSOR^{ACG}_{pl}$ & 0.305 & 0.730 & 0.233 & 0.212 & 0.191 \\
$SSOR^{ACG}_{nl}$ & \textbf{\textit{0.147}} & \textbf{\textit{0.086}} & \textbf{\textit{0.057}} & 0.034 & 0.021 \\
\hline
\end{tabular}
}
}
\qquad\qquad
\subtable[$\overline{Correct}^{\%}_{Method}$\label{PercentageSensitivityComparisonsUnique}]{%
\resizebox{0.43\textwidth}{!}{
\begin{tabular}{lccccc}
\hline
Model            & $z=4$ & $z=9$ & $z=14$ & $z=19$ & $z=25$    \\   
\hline 
$Uniform$ & 67.98\% & 59.44\% & 40.89\% & 12.31\% & 1.73\% \\
$Logistic$ & 61.06\% & 76.79\% & 84.94\% & \textbf{\textit{90.67\%}} & \textbf{\textit{94.37\%}} \\
$SSOR^{ACG}_{pl}$ & 78.48\% & 81.18\% & 78.48\% & 70.33\% & 38.30\% \\
$SSOR^{ACG}_{nl}$ & \textbf{\textit{85.23\%}} & \textbf{\textit{89.86\%}} & \textbf{\textit{89.87\%}} & 87.01\% & 91.23\% \\
\hline
\end{tabular}
}
}
\end{center}
\end{table}

\begin{table}[!h]
\begin{center}
\caption{Mean distance between the PWI matrices ($\overline{d}(PWI_{DM},PWI_{Method})$ and mean percentage of correct comparisons of no reference pairs of alternatives ($\overline{Correct}^{\%}_{Method}$) considering the distribution ${\mathbf{p}}^{DM}=\mathbf{p}^{Nor}$}\label{SensitivityComparisonsNorm}
\subtable[$\overline{d}(PWI_{DM},PWI_{Method})$\label{DistanceSensitivityComparisonsNorm}]{%
\resizebox{0.4\textwidth}{!}{
\begin{tabular}{lccccc}
\hline
Model            & $z=4$ & $z=9$ & $z=14$ & $z=19$ & $z=25$    \\     
\hline
$Uniform$ & 0.226 & 0.236 & 0.231 & 0.237 & 0.233 \\
$SSOR^{ACG}_{pl}$ & 0.181 & 0.149 & 0.127 & 0.118 & 0.114 \\
$SSOR^{ACG}_{nl}$ & \textbf{\textit{0.118}} & \textbf{\textit{0.101}} & \textbf{\textit{0.095}} & \textbf{\textit{0.095}} & \textbf{\textit{0.094}} \\
\hline
\end{tabular}
}
}
\qquad\qquad
\subtable[$\overline{Correct}^{\%}_{Method}$\label{PercentageSensitivityComparisonsNorm}]{%
\resizebox{0.46\textwidth}{!}{
\begin{tabular}{lccccc}
\hline
Model            & $z=4$ & $z=9$ & $z=14$ & $z=19$ & $z=25$    \\   
\hline 
$Uniform$  & 77.44\% & 70.05\% & 60.44\% & 37.80\% & 7.87\% \\
$SSOR^{ACG}_{pl}$ & 86.85\% & 89.55\% & 89.79\% & 83.19\% & 60.50\% \\
$SSOR^{ACG}_{nl}$ & \textbf{\textit{89.13\%}} & \textbf{\textit{91.86\%}} & \textbf{\textit{94.51\%}} & \textbf{\textit{94.38\%}} & \textbf{\textit{93.23\%}} \\
\hline
\end{tabular}
}
}
\end{center}
\end{table}
Looking at the data in Tables \ref{SensitivityComparisonsUnique} and \ref{SensitivityComparisonsNorm} the following can be observed: 
\begin{itemize}
\item Considering $\delta_{U_{ref}}$ as the probability distribution of the artificial DM over $\mathbf{\Omega}$, we can observe that $SSOR^{ACG}_{nl}$ is the best among the considered methods if the number of pairwise comparisons provided by the DM is 4, 9 or 14. This is supported by the data in Tables \ref{DistanceSensitivityComparisonsUnique} and \ref{PercentageSensitivityComparisonsUnique} reporting the mean distance between the PWI matrices and the mean percentage of correct comparisons of no reference pairs of alternatives, respectively. In the other two cases ($z=19$ and $z=25$) the $Logistic$ method is the best both considering the mean distance and the mean percentage even if the difference with the equivalent data reported for $SSOR^{ACG}_{nl}$ is not very large. \\
An interesting thing to observe is that increasing the number of pairwise comparisons of reference alternatives provided by the artificial DM is not always beneficial for $SSOR^{ACG}_{pl}$. Indeed, passing from $z=4$ to $z=9$, the mean percentage of correct comparisons of no reference pairs of alternatives increases, while, it drastically decreases passing from $z=19$ to $z=25$ of both distributions. This is mainly due to the fact that, as explained in Section \ref{goff}, the probability distribution built by $SSOR^{ACG}_{pl}$ is a piecewise linear function defined by three breakpoints only. Therefore, since the number of variables defining the probability distribution is quite small, increasing the amount of preference information provided by the DM increases the infeasibility of the LP problem \ref{lpp2}\footnote{If $\mathbf{p}^{DM}=\delta_{U_{ref}}$, the percentage of infeasible LP problems (\ref{lpp2}) considering $z=19$ and $z=25$ is $56.10\%$ and $71.90\%$, respectively.}. Therefore, the probability distribution obtained solving the LP problem is such that not only the comparisons between no reference pairs of alternatives are restored, but also some of the pairwise comparisons of the reference alternatives given by the artificial DM. \\
The $Uniform$ method is the worst among the considered methods for all considered values of $z$;
\item assuming that the DM's probability distribution over $\mathbf{\Omega}$ is $\mathbf{p}^{Nor}$, we can observe from Tables \ref{DistanceSensitivityComparisonsNorm} and \ref{PercentageSensitivityComparisonsNorm} that $SSOR^{ACG}_{nl}$ is the best method for all considered values of $z$. In particular, we can observe that the increase in the number of pairwise comparisons provided by the DM does involve always an increase in the mean percentage of correct comparisons of the no reference pairs of alternatives even if, in these cases (passing from $z=14$ to $z=19$ and, then, to $z=25$) the diminishing on the percentage is quite small. 
Also in this case, increasing the amount of preference information is not beneficial for the $SSOR^{ACG}_{pl}$ method since for $z=25$, it presents the lowest percentage of correct comparisons of no reference pairs of alternatives\footnote{If $\mathbf{p}^{DM}=\mathbf{p}^{Nor}$, the percentage of infeasible LP problems (\ref{lpp2}) considering $z=19$ and $z=25$ is $40.30\%$ and $50.30\%$, respectively.}.
\end{itemize}

%%%%%%%%%%%%%%%%%%%%%%%%%%%%%%%%%%%%%%%
\section{Discussion}\label{Comments}%%%
%%%%%%%%%%%%%%%%%%%%%%%%%%%%%%%%%%%%%%%
The large amount of experiments we performed show that, undoubtedly, $SSOR^{ACG}_{nl}$ is the best among the methods taken into account and, therefore, it should be applied in practice. In most of the problems considered varying the number of alternatives or criteria, as well as the number of pairwise comparisons provided by the DM or the assumed probability distribution on $\mathbf{\Omega}$, $SSOR^{ACG}_{nl}$ presents the best results in terms of distance from the PWI or RAI matrices and in terms of percentage of correct comparisons of no reference pairs of alternatives. For this reason, we shall comment, in the following the $SSOR^{ACG}_{nl}$ method only. 

Regarding the choice of the barycenter as the reference model, we have proved that a different choice negatively affects the obtained results. In particular, the two versions of the $SSOR^{ACG}_{nl}$ method in which the reference model is unknown presented almost always the worst results. This probably has to be interpreted in terms of overfitting due to a larger number of variables to be estimated on the basis of the available data. Moreover, with respect to the sensitivity of the results to the choice of the barycenter as the reference model, we have also proved that by replacing the barycenter with another value function in its neighborhood, the results do not change in a considerable way. This proves that the choice of the barycenter as reference model is also robust with respect to the recommendations obtained by the method. 

An aspect that should be underlined is that the goodness of the results is not strictly dependent on the amount of preference information provided by the DM. Indeed, as proved in Section \ref{Sensitivityz}, on the one hand, the increase in the number of pairwise comparisons involves, in general, a diminishing of the distance from the PWI matrices, while, on the other hand, this is not the case for the percentage of correct comparisons of no reference pairs of alternatives. Indeed, passing from $z=4$ to $z=9$ and, then, to $z=14$, we observed an increase in the percentage of correct comparisons of no reference pairs of alternatives but, then, passing from $z=14$ to $z=19$ we got a deterioration of the same percentage that continued also passing to $z=25$ in considering one of the two assumed probability distributions. This proves that research has to be devoted to checking the ``optimal" amount of preference information that should be required to the DM in this case or how to modify the proposed algorithms when the amount of preference information increases.

%%%%%%%%%%%%%%%%%%%%%%%%%%%%%%%%%%
\section{Conclusions}\label{conc}%
%%%%%%%%%%%%%%%%%%%%%%%%%%%%%%%%%%
In this paper, we propose some methodologies to infer a probability distribution defined over a set of models compatible with preference information provided by the Decision Maker (DM) and sampled from the corresponding space. Differently from the Subjective Stochastic Ordinal Regression proposed by \cite{CorrenteEtAl2016KNOSYS}, in our method the mass attached to each compatible model is dependent on the distance from a reference model representing the basic preference tendency of the DM to evaluate the alternatives under consideration. In particular, the mass function is a non-increasing function of the distance from the reference model so that, the closest the compatible model to the reference one, the greatest the mass attached to it. To get the mass function in a first approach considering only the constraints ensuring that the larger the distance from the reference model the smaller the probability mass, a simple LP problem or the optimization of a single variable function needs to be solved. \\
To make the computation of the mass function simpler, in a second proposal, we assume that it is a piecewise linear function of the same distance from the reference model. Considering a few reference distance values and associating them with a mass, the mass attached to all the other sampled models is obtained by linear interpolation. In this way, the considered mass function is defined only by the mass of the reference distances that can be computed again by solving a simple LP problem. \\
Finally, in our third proposal, we assume that the probability distribution on the sample of compatible models takes the normal or the exponential form. In this case, the computation of the considered probability distribution is based on the solution of a non-linear programming problem optimizing a single variable function. In the three proposals, we assumed that the reference model used to compute the probability mass of each function is the barycenter of the space of sampled value functions compatible with some preferences given by the DM.

To prove the effectiveness and reliability of the three proposals, we performed an extensive set of simulations assuming the existence of an artificial DM whose probability distribution on the sample of value functions needs to be estimated. For such a reason, different forms of DM's probability distributions have been taken into account. For each of the considered methods, we computed the distance between rank acceptability indices and pairwise winning indices matrices computed using the artificial DM's probability distribution and the same matrices computed considering the estimated probability distribution. To check if the distance values are significant from the statistical point of view we performed two versions of the 2-sample Kolmogorov-Smirnov test. 

In addition to these experiments, we made also an extensive sensitivity analysis with respect to several elements affecting the results obtained by the application of the methodology we are proposing.  In particular, we tested the quality of the results provided by our approach to the variation of the number of alternatives or criteria, the number of pairwise comparisons provided by the DM, and the choice of the barycenter as the reference model. Since we proved that the barycenter is the optimal choice for the reference model, we also performed a sensitivity analysis on it showing that choosing another value function in its neighborhood as reference model does not affect the goodness of the obtained results. This proves the robustness of this choice.  

Both the results and the tests show that the third proposal assuming the existence of a probability distribution of a normal or exponential form on the sample of value functions compatible with the preferences given by the DM is the best among those under consideration. For such a reason, we suggest using it in real-world applications. 

On the basis of the methodology proposed in this paper, we envisage the following future developments. First of all, statistical properties of the distribution inferred by the proposed models should be studied. Secondly, it is worthwhile to apply similar approaches to preference models different from the weighted sum mainly considered in this paper. This is the case of the piecewise linear additive value functions, the non-additive multicriteria aggregation procedures such as the Choquet integral, the outranking methods such as ELECTRE \citep{FigueiraEtAl2013,FigueiraMousseauRoy2016} or PROMETHEE methods \citep{BransVincke1985}. After, we believe that it is relevant to investigate the possibility of improving the obtained results with an interactive procedure updating the probability in the space of feasible models on the basis of preferences iteratively elicited from the DM during the decision-aiding procedure. Moreover, one can imagine applying the proposed procedure in interactive multiobjective optimization methods \citep{BrankeEtAl2008}, both evolutionary or not, in order to guide the algorithms to discover the most preferred solutions in the Pareto front.

%%%%%%%%%%%%%%%%%%%%%%%%%%%%
\section*{Acknowledgments}%%
%%%%%%%%%%%%%%%%%%%%%%%%%%%%
The authors wish to acknowledge the support of the Ministero dell'Istruzione, dell'Universit\'{a} e della Ricerca (MIUR) - PRIN 2017, project ``Multiple Criteria Decision Analysis and Multiple Criteria Decision Theory'', grant 2017CY2NCA. Moreover, the authors wish to acknowledge Prof. Antonio Punzo for his valuable support on the performed statistical tests.  

%%%%%%%%%%%%%%%%%%%%%%%%%%%%%%%%%%%%%%%%%%%%%%%%%%%%%%%%%%
\bibliographystyle{plainnat}                             %
%\bibliography{Full_bibliography,Full_bibliographySort}  %
\bibliography{BiBliography_Final}                              %
%\vspace{0,4cm}                                          %
%%%%%%%%%%%%%%%%%%%%%%%%%%%%%%%%%%%%%%%%%%%%%%%%%%%%%%%%%%
%\noindent\textbf{Proof of Proposition }}\\%%
%%%%%%%%%%%%%%%%%%%%%%%%%%%%%%%%%%%%%%%%%%%%%%%%%%%%%%%%%%

\end{document}